\documentclass[10pt,a4paper,english,french,british]{article}
\usepackage[T1]{fontenc}
\usepackage[utf8x]{inputenc}
\usepackage{amsmath}
\usepackage{graphicx}
\usepackage{amssymb}
\usepackage{esint}
\usepackage[authoryear]{natbib}

\makeatletter

\providecommand{\tabularnewline}{\\}

\usepackage{longtable}\usepackage{nicefrac}\usepackage{bm}\usepackage{amsfonts}\usepackage{mathrsfs}\usepackage{gensymb}\usepackage{amsthm}\usepackage{xcolor}\usepackage{bbm}
\usepackage{xspace}

\newtheorem{thme}{Theorem}
\newtheorem{assumption}{A}
\newtheorem{lem}[thme]{Lemma}
\newtheorem{prop}[thme]{Proposition}
\newtheorem{req}[thme]{Remark}
\newtheorem{res}[thme]{Result}
\newtheorem*{berb}{Berbee's coupling lemma}
\newtheorem*{talagrand}{Talagrand's inequality}
\newtheorem*{burkholder}{Burkholder Davis Gundy inequality}

\DeclareMathOperator{\variance}{Var}
\newcommand{\norm}[1]{\Vert #1 \Vert}
\newcommand{\rond}[1]{\mathscr{#1}}
\newcommand{\E}[1]{\mathbb{E}\left(#1\right)}
\newcommand{\EE}[1]{\mathbb{E}\left[#1\right]}
\newcommand{\var}[1]{\variance\left(#1\right)}

\newcommand{\units}[1]{\mathbbmss{1}_{#1}}

\newcommand{\burk}{Burkholder Davis Gundy inequality\xspace}
\title{Non parametric estimation of the diffusion coefficents of a diffusion with jumps}
\author{Émeline Schmisser}
\date{}

\usepackage{babel}\addto\extrasfrench{%
   \providecommand{\fg}{\ifdim\lastskip>\z@\unskip\fi~\frqq}%
}

\usepackage{babel}\addto\extrasfrench{%
   \providecommand{\fg}{\ifdim\lastskip>\z@\unskip\fi~\frqq}%
}

\makeatother

\usepackage{babel}\addto\extrasfrench{%
   \providecommand{\fg}{\ifdim\lastskip>\z@\unskip\fi~\frqq}%
}

\makeatother

\usepackage{babel}
\addto\extrasfrench{\providecommand{\fg}{\ifdim\lastskip>\z@\unskip\fi~\frqq}}

\begin{document}

\title{Non parametric estimation of the coefficients of a diffusion with
jumps}
\maketitle
\begin{abstract}
In this article, we consider a jump diffusion process $\left(X_{t}\right)_{t\geq0}$,
with drift function $b$, diffusion coefficient $\sigma$ and jump
coefficient $\xi^{2}$. This process is observed at discrete times
$t=0,\Delta,\ldots,n\Delta$. The sampling interval $\Delta$ tends
to 0 and $n\Delta$ tends to infinity. We assume that $\left(X_{t}\right)_{t\geq0}$
is ergodic, strictly stationary and exponentially $\beta$-mixing.
We use a penalized least-square approach to compute adaptive estimators
of the functions $\sigma^{2}+\xi^{2}$ and $\sigma^{2}$. We provide
bounds for the risks of the two estimators. 
\end{abstract}

\selectlanguage{french}%
\begin{abstract}
Nous observons une diffusion à sauts $(X_{t})_{t\geq0}$ à des instants
discrets $t=0,\Delta,\ldots,n\Delta$. Le temps d'observation $n\Delta$
tend vers l'infini et le pas d'observation $\Delta$ tend vers 0).
Nous supposons que le processus $(X_{t})_{t\geq0}$ est ergodique,
stationnaire et exponentiellement $\beta$-mélangeant. Nous construisons
des estimateurs adaptatifs des fonctions $\sigma^{2}+\xi^{2}$ et
$\sigma^{2}$, où $\sigma^{2}$ est le coefficient de diffusions et
$\xi^{2}$ le coefficient de sauts, grâce à une méthode de moindres
carrés pénalisés. Nous majorons le risque de ces estimateurs de manière
non asymptotique.
\end{abstract}
\selectlanguage{british}%

\section{Introduction}

We consider the stochastic differential equation (SDE): \begin{equation}
dX_{t}=b(X_{t^{-}})dt+\sigma(X_{t^{-}})dW_{t}+\xi(X_{t^{-}})dL_{t},\qquad X_{0}=\eta\label{eq:EDS}\end{equation}
 with $\eta$ a random variable, $(W_{t})_{t\geq0}$
a Brownian motion independent of $\eta$ and $\left(L_{t}\right)_{t\geq0}$
a pure jump centered Lévy process independent of $\left(\left(W_{t}\right)_{t\geq0},\eta\right)$:

\[
L_{t}=\int_{0}^{t}\int_{z\in\mathbb{R}}z\left(\mu(dt,dz)-\nu(dz)\right)dt\]
 where $\mu$ is a Poisson measure of intensity $\nu(dz)dt$, with
$\int_{\mathbb{R}}(z^{2}\wedge1)\nu(dz)<\infty$. The process $(X_{t})_{t\geq0}$
is assumed to be ergodic, stationary and exponentially $\beta$-mixing.
It is observed at discrete times $t=0,\Delta,\ldots,n\Delta$ where
the sampling interval $\Delta$ tends to 0 and the time of observation
$n\Delta$ tends to infinity. Our aim is to construct adaptive non-parametric
estimators of $\xi^{2}$ and $\sigma^{2}$ on a compact set $A$.

Diffusions with jumps become powerful tools to model processes in
biology, physics, social sciences, medical sciences, economics, and
a variety of financial applications such as interest rate modelling
or derivative pricing. However, if the non-parametric estimation of
the coefficients of a diffusion without jumps is well known (see for
instance \citet{hofmann99} or \citet{comtegenon2007}), to our knowledge,
there do not exist adaptive estimators for the coefficients of a jump
diffusion, neither minimax rates of convergence. \citet{shimizu2008}
construct maximum-likelihood parametric estimators of $\sigma^{2}$
and $\xi^{2}$. Their estimators converge with rates $\sqrt{n}$ and
$\sqrt{n\Delta}$ respectively. \citet{mancinireno2011} and \citet{muhammad_wang_lin}
construct non-parametric estimators of $\sigma^{2}$ and $\sigma^{2}+\xi^{2}$
thanks to kernel or local polynomials estimators. The estimator of
$\sigma^{2}$ converges with rate $\sqrt{hn}$, meanwhile the estimator
of $\xi^{2}+\sigma^{2}$ converges with rate $\sqrt{n\Delta h}$,
where $h$ is the bandwidth of the estimator.

In this paper, we construct non-parametric estimators of $g=\sigma^{2}+\xi^{2}$
and $\sigma^{2}$ under the asymptotic framework $n\Delta\rightarrow\infty$
and $\Delta\rightarrow0$ by model selection. This method was introduced
by \citet{birgemassart98}. We consider first the following random
variables \[
T_{k\Delta}=\frac{(X_{(k+1)\Delta}-X_{k\Delta})^{2}}{\Delta}=\sigma^{2}(X_{k\Delta})+\xi^{2}(X_{k\Delta})+\textrm{noise}+\textrm{remainder.}\]
We introduce a sequence of increasing subspaces $S_{m}$ of $L^{2}(A)$
and we construct a sequence of estimators $\hat{g}_{m}$ by minimizing
over each $S_{m}$ a contrast function \[
\gamma_{n}(t)=\frac{1}{n}\sum_{k=1}^{n}(T_{k\Delta}-t(X_{k\Delta}))^{2}.\]
We bound the risk of $\hat{g}_{m}$, then we introduce a penalty function
$pen(m)$ and me minimize on $m$ the function $\gamma_{n}(\hat{g}_{m})+pen(m)$.
If the Lévy measure $\nu$ is sub-exponential, the adaptive estimator
$\hat{g}_{\hat{m}}$ satisfies an oracle inequality (up to a multiplicative
constant). 

To estimate the function $\sigma^{2}$, we need to cut off the jumps.
We minimize over each $S_{m}$ the contrast function \[
\tilde{\gamma}_{n}(t)=\frac{1}{n}\sum_{k=1}^{n}\left(T_{k\Delta}\units{\left|X_{(k+1)\Delta}-X_{k\Delta}\right|\leq C_{\Delta}}-t(X_{k\Delta})\right)^{2}\quad\textrm{where}\quad C_{\Delta}\propto\sqrt{\Delta}\ln(n).\]
We obtain a sequence of estimators $\hat{\sigma}_{m}^{2}$ of $\sigma^{2}$.
The risk of these estimators depends on the Blumenthal-Getoor index
of $\nu$. To construct an adaptive estimator, $\hat{\sigma}_{\hat{m}}^{2}$,
we again introduce a penalty function $\widetilde{pen}(m)$. The estimator
$\hat{\sigma}_{\hat{m}}^{2}$ automatically realizes a bias-variance
compromise. The rates of convergence obtained for $\hat{g}_{m}$ and
$\hat{\sigma}_{m}^{2}$ are similar to those obtained by \citet{muhammad_wang_lin}
and \citet{mancinireno2011}. 

This article is composed as follows:~in Section \ref{sec:Model},
we specify the model and its assumptions. In Sections \ref{sec:Estimation_sigma_xi}
and \ref{sec:Estimation-sigma}, we construct the estimators and bound
their risks. Section \ref{sec:Simulations} is devoted to the simulations
and proofs are gathered in Section \ref{sec:Proofs}.

\section{Model \label{sec:Model}}

We consider the stochastic differential equation \eqref{eq:EDS}.
We assume that the following assumptions are fulfilled:

\begin{assumption} \label{hypo_modele} 
\begin{enumerate}
\item \label{enu:hypo_lip}The functions $b$, $\sigma$ and $\xi$ are
Lipschitz. 
\item \label{enu:hypo_majoration_sigma_xi}The functions $\sigma$ and $\xi$
are bounded: $\exists\sigma_{0}^{2},\xi_{0}^{2}$ such that \[
\forall x\in\mathbb{R},\:0<\sigma^{2}(x)\leq\sigma_{0}^{2}\quad\textrm{and}\quad0<\xi^{2}(x)\leq\xi_{0}^{2}.\]
 Moreover either there exists a positive constant constant $\sigma_{1}^{2}$
such that $\forall x\in\mathbb{R},\:\sigma^{2}(x)\geq\sigma_{1}^{2}>0$,
or there exists $\xi_{1}^{2}$ such that, $\forall x\in\mathbb{R},\:\xi^{2}(x)\geq\xi_{1}^{2}>0$. 
\item \label{enu:hypo_b_elastique}The function $b$ is elastic: $\exists M,C$,
$\forall x,\vert x\vert>M$, $b(x)x\leq-Cx^{2}$. 
\item \label{enu:hypo_blumenthal}The Lévy measure $\nu$ satisfies: \[
\nu\left(\left\{ 0\right\} \right)=0\quad\textrm{and}\quad\int_{-\infty}^{+\infty}z^{8}\nu(dz)<\infty\]
 and the Blumenthal-Getoor index is strictly less than 2: there exists
$\beta\in[0,2[$ such that $\int_{-1}^{1}z^{\beta}\nu(dz)<\infty$.
This is a classical assumption (see for instance \citet{mai2012}).
In order to ensure the uniqueness of the function $\xi$, we also
assume that $\int_{-\infty}^{+\infty}z^{2}\nu(dz)=1$. 
\end{enumerate}
\end{assumption}

If Assumption A\ref{hypo_modele}.\ref{enu:hypo_lip} is satisfied,
SDE \eqref{eq:EDS} as a unique solution. According to \citet{masuda2007},
under assumptions A\ref{hypo_modele}.(\ref{enu:hypo_lip}-\ref{enu:hypo_b_elastique}),
the process $\left(X_{t}\right)_{t\geq0}$ is exponentially $\beta$-mixing
and has a unique invariant probability. Moreover, under assumption
A\ref{hypo_modele}.(\ref{enu:hypo_blumenthal}), $\E{X_{t}^{8}}<\infty$.
Then we can assume:~

\begin{assumption}\label{hypo_stationnaire.}

\label{enu:hypo_stationnarite}The process $\left(X_{t}\right)_{t\geq0}$
is stationary, exponentially $\beta-$mixing and its stationary measure
has a density $\pi$ which is bounded on any compact set. 

\end{assumption}

The following result is very useful. It comes from \citet{dellacheriemeyer}
or \citet{applebaum}.

\begin{burkholder}

Let us consider the filtration \[
\rond{F}_{t}=\sigma(\eta,(W_{s})_{0\leq s\leq t},(L_{s})_{0\leq s\leq t}).\]
 Then, for any $p\geq2$, there exists a constant $C_{p}>0$ such
that: \[
\E{\sup_{s\in[t,t+h]}\left|\int_{s}^{t}\sigma(X_{u})dW_{u}\right|^{p}\left|\rond{F}_{t}\right.}\leq C_{p}\E{\left|\int_{t}^{t+h}\sigma^{2}(X_{u})du\right|^{p/2}\left|\rond{F}_{t}\right.}\]
 and \begin{align*}
\E{\sup_{s\in[t,t+h]}\left|\int_{s}^{t}\xi(X_{u^{-}})dL_{u}\right|^{p}\left|\rond{F}_{t}\right.} & \leq C_{p}\E{\left|\int_{t}^{t+h}\xi^{2}(X_{u})du\right|^{p/2}\left|\rond{F}_{t}\right.}\\
 & +C_{p}\E{\left|\int_{t}^{t+h}\left|\xi^{p}(X_{u})\right|du\right|\left|\rond{F}_{t}\right.}\int_{\mathbb{R}}\left|z\right|^{p}\nu(dz)\end{align*}

\end{burkholder}

The following proposition derives from this result.

\begin{prop}\label{prop_Xk}

For any integer $p$ and any $t\leq1$:~

\[
\E{\sup_{0\leq s\leq t}\left(X_{t+u}-X_{u}\right)^{2p}}\lesssim t.\]

\end{prop}

Now we introduce an increasing sequence of vectorial subspaces $\left(S_{m}\right)_{m\geq0}$
of $L^{2}(A)$ satisfying the following properties:~

\begin{assumption}\label{hypo_espaces} 
\begin{enumerate}
\item The subspaces $S_{m}$ have finite dimension $D_{m}$ and are increasing:~
$\forall m$, $S_{m}\subseteq S_{m+1}$. 
\item The $\norm{.}_{L^{2}}$ and $\norm{.}_{\infty}$ norms are connected:~
\[
\exists\phi_{1},\forall m,\forall t\in S_{m},\quad\norm{t}_{\infty}^{2}\leq\phi_{1}D_{m}\norm{t}_{L^{2}}^{2}\]
 with $\norm{t}_{L^{2}}^{2}=\int_{A}t^{2}(x)dx$ and $\norm{t}_{\infty}=\sup_{x\in A}\left|t(x)\right|$. 
\item For any function $t\in\rond{B}_{2,\infty}^{\alpha}$, \[
\exists c,\:\forall m,\:\norm{t-t_{m}}_{L^{2}}^{2}\leq cD_{m}^{-2\alpha}\]
 where $t$ is the orthogonal projection $L^{2}$ of $t$ on $S_{m}$. 
\end{enumerate}
\end{assumption}

The vectorial subspaces generated by the trigonometric polynomials,
the piecewise polynomials, the spline functions and the wavelets satisfy
these properties (see \citet{meyer} and \citet{devorelorentz} for
the proofs).

\section{Estimation of $\sigma^{2}+\xi^{2}.$ \label{sec:Estimation_sigma_xi}}

Let us set $Z_{k\Delta}=\int_{k\Delta}^{(k+1)\Delta}\sigma(X_{s})dW_{s}$
and $J_{k\Delta}=\int_{k\Delta}^{(k+1)\Delta}\xi(X_{s^{-}})dL_{s}$.
To estimate $\sigma^{2}$ for a diffusion process (without jumps),
we can consider the random variables \[
T_{k\Delta}=\frac{(X_{(k+1)\Delta}-X_{k\Delta})^{2}}{\Delta}\]
 (see \citet{comtegenon2007}). For jump diffusions, \[
X_{(k+1)\Delta}-X_{k\Delta}=\int_{k\Delta}^{(k+1)\Delta}b(X_{s})ds+Z_{k\delta}+J_{k\Delta}\]
 and therefore

\[
T_{k\Delta}=\sigma^{2}(X_{k\Delta})+\xi^{2}(X_{k\Delta})+A_{k\Delta}+B_{k\Delta}+E_{k\Delta}\]
 where \begin{align*}
A_{k\Delta} & =A_{k\Delta}^{(1)}+A_{k\Delta}^{(2)}+A_{k\Delta}^{(3)}+A_{k\Delta}^{(4)}\\
 & =\frac{1}{\Delta}\left(\int_{k\Delta}^{(k+1)\Delta}b(X_{s})ds\right)^{2}+\frac{2}{\Delta}\left(Z_{k\Delta}+J_{k\Delta}\right)\int_{k\Delta}^{(k+1)\Delta}\left(b(X_{s})-b(X_{k\Delta})\right)ds\\
 & +\frac{1}{\Delta}\int_{k\Delta}^{(k+1)\Delta}\left(\sigma^{2}(X_{s})-\sigma^{2}(X_{k\Delta})\right)ds+\frac{1}{\Delta}\int_{k\Delta}^{(k+1)\Delta}\left(\xi^{2}(X_{s})-\xi^{2}(X_{k\Delta})\right)ds,\end{align*}
 \[
B_{k\Delta}=B_{k\Delta}^{(1)}+B_{k\Delta}^{(2)}=2b(X_{k\Delta})Z_{k\Delta}+\frac{1}{\Delta}\left[Z_{k\Delta}^{2}-\int_{k\Delta}^{(k+1)\Delta}\sigma^{2}(X_{s})ds\right]\]
 and \begin{align*}
E_{k\Delta} & =E_{k\Delta}^{(1)}+E_{k\Delta}^{(2)}+E_{k\Delta}^{(3)}\\
 & =2b(X_{k\Delta})J_{k\Delta}+\frac{2}{\Delta}Z_{k\Delta}J_{k\Delta}+\frac{1}{\Delta}\left[J_{k\Delta}^{2}-\int_{k\Delta}^{(k+1)\Delta}\xi^{2}(X_{s^{-}})ds\right].\end{align*}
 The term $A_{k\Delta}$ is small, whereas $B_{k\Delta}$ and $E_{k\Delta}$
are centred. The random variables $B_{k\Delta}$ depend on the Brownian
motion $\left(W_{t}\right)_{t\geq0}$, while $E_{k\Delta}$ depends
on the jump process $\left(L_{t}\right)_{t\geq0}$. The following
lemma is derived from Proposition \ref{prop_Xk} and the \burk.

\begin{lem}\label{lem_majoration_AB} 
\begin{itemize}
\item $\E{A_{k\Delta}^{2}}\lesssim\Delta$ and $\E{A_{k\Delta}^{4}}\lesssim\Delta.$ 
\item $\E{B_{k\Delta}\left|\rond{F}_{k\Delta}\right.}=0$, $\E{B_{k\Delta}^{2}\left|\rond{F}_{k\Delta}\right.}\lesssim1$
and $\E{B_{k\Delta}^{4}}\lesssim1$. 
\item $\E{E_{k\Delta}\left|\rond{F}_{k\Delta}\right.}=0$, $\E{E_{k\Delta}^{2}\left|\rond{F}_{k\Delta}\right.}\lesssim1/\Delta$
and $\E{E_{k\Delta}^{4}}\lesssim1/\Delta^{3}$. 
\end{itemize}
\end{lem}

\subsection{Estimation for fixed $m$}

For any $m\in\rond{M}_{n}=\left\{ m,\: D_{m}\leq\rond{D}_{n}\right\} $
where the maximal dimension $\rond{D}_{n}$ satisfies $\rond{D}_{n}\leq\sqrt{n\Delta}/\ln(n)$,
we construct an estimator $\hat{g}_{m}$ of $g=\sigma^{2}+\xi^{2}$
by minimizing on $S_{m}$ the contrast function \[
\gamma_{n}(t)=\frac{1}{n}\sum_{k=1}^{n}\left(t(X_{k\Delta})-T_{k\Delta}\right)^{2}.\]
 Let us bound the empirical risk $\rond{R}_{n}(\hat{g}_{m})$, where
\[
\rond{R}_{n}(t)=\E{\norm{t-g}_{n}^{2}}\quad\textrm{with}\quad\norm{t}_{n}^{2}=\frac{1}{n}\sum_{k=1}^{n}t^{2}(X_{k\Delta}).\]
 We set $\norm{t}_{\pi}^{2}=\int_{A}t^{2}(x)\pi(x)dx$ and $g_{A}=g\units{A}$. 

We have that \begin{align*}
\gamma_{n}(t) & =\frac{1}{n}\sum_{k=1}^{n}\left(t(X_{k\Delta})-g(X_{k\Delta})+A_{k\Delta}+B_{k\Delta}+E_{k\Delta}\right)^{2}\\
 & =\norm{t-g}_{n}^{2}+\frac{1}{n}\sum_{k=1}^{n}(A_{k\Delta}+B_{k\Delta}+E_{k\Delta})^{2}\\
 & -\frac{2}{n}\sum_{k=1}^{n}\left(A_{k\Delta}+B_{k\Delta}+E_{k\Delta}\right)\left(g(X_{k\Delta})-t(X_{k\Delta})\right).\end{align*}
 As $\hat{g}_{m}$ minimizes $\gamma_{n}(t)$, the inequality $\gamma_{n}(\hat{g}_{m})\leq\gamma_{n}(g_{m})$
holds and then \[
\norm{\hat{g}_{m}-g}_{n}^{2}\leq\norm{g_{m}-g}_{n}^{2}+\frac{2}{n}\sum_{k=1}^{n}\left(A_{k\Delta}+B_{k\Delta}+E_{k\Delta}\right)\left(\hat{g}_{m}(X_{k\Delta})-g_{m}(X_{k\Delta})\right).\]
 By Cauchy-Schwarz, and as $\hat{g}_{m}$ and $g_{m}$ are $A$-supported,
\[
\norm{\hat{g}_{m}-g_{A}}_{n}^{2}\leq\norm{g_{m}-g_{A}}_{n}^{2}+\frac{12}{n}\sum_{k=1}^{n}A_{k\Delta}^{2}+\frac{1}{12}\norm{\hat{g}_{m}-g_{m}}_{n}^{2}+12\sup_{t\in\rond{B}_{m}}\nu_{n}^{2}(t)+\frac{1}{12}\norm{\hat{g}_{m}-g_{m}}_{\pi}^{2}\]
 where $\rond{B}_{m}=\left\{ t\in S_{m},\norm{t}_{\pi}^{2}\leq1\right\} $
and $\nu_{n}(t)=\frac{1}{n}\sum_{k=1}^{n}(B_{k\Delta}+E_{k\Delta})t(X_{k\Delta})$.
Let us set \[
\Omega_{n}=\left\{ \omega,\:\forall m\in\rond{M}_{n},\forall t\in S_{m},\left|\frac{\norm{t}_{n}^{2}}{\norm{t}_{\pi}^{2}}-1\right|\leq\frac{1}{2}\right\} \]
 where the norms $\left\Vert .\right\Vert _{\pi}$ and $\left\Vert .\right\Vert _{n}$
are equivalent. The following lemma is proved by \citet{comtegenon2007}
for diffusion processes, but only relies of the $\beta$-mixing and
stationary properties. 

\begin{lem}\label{majoration_omegaC}

\[
\mathbb{P}\left(\Omega_{n}^{c}\right)\leq\frac{c}{n^{8}}.\]

\end{lem}

We obtain that

\[
\E{\norm{\hat{g}_{m}-g_{A}}_{n}^{2}\units{\Omega_{n}}}\leq3\norm{g_{m}-g_{A}}_{\pi}^{2}+12\E{A_{k\Delta}^{2}}+12\E{\sup_{t\in\rond{B}_{m}}\nu_{n}^{2}(t)}.\]

On $\Omega_{n}$, any function $t\in S_{m}$ satisfies: $\norm{t}_{\pi}^{2}\leq2\norm{t}_{n}^{2}$.
Moreover, for any deterministic function $t$, $\E{\norm{t}_{n}^{2}}=\norm{t}_{\pi}^{2}$.
Consequently:

\[
\E{\norm{\hat{g}_{m}-g_{A}}_{n}^{2}\units{\Omega_{n}}}\leq3\norm{g_{m}-g_{A}}_{\pi}^{2}+12\E{A_{\Delta}^{2}}+12\E{\sup_{t\in\rond{B}_{m}}\nu_{n}^{2}(t)}.\]
By Assumption A\ref{hypo_stationnaire.}, $\pi$ is bounded on $A$
and then $\norm{g_{m}-g_{A}}_{\pi}^{2}\lesssim\norm{g_{m}-g_{A}}_{L^{2}}^{2}.$
The remainder of the proof is done in Section \ref{sec:Proofs}. 

\begin{thme}\label{thme_risque_est_sigma_xi_m_fixe}

Under Assumptions A\ref{hypo_modele}-A\ref{hypo_espaces}, if $m\in\rond{M}_{n}$,
the risk of the estimator $\hat{g}_{m}$ is bounded by: \[
\rond{R}_{n}(\hat{g}_{m})\lesssim\norm{g_{m}-g_{A}}_{L^{2}}^{2}+\frac{D_{m}}{n\Delta}\xi_{0}^{4}+\frac{D_{m}}{n}(\sigma_{0}^{4}+\sigma_{0}^{2}\xi_{0}^{2})+\frac{1}{n\Delta}+\Delta\]
 where $g_{m}$ is the orthogonal $(L^{2})$ projection of $g$ on
$S_{m}$.

\end{thme}

We have to find a good compromise between the bias term, $\norm{g_{m}-g_{A}}_{L^{2}}^{2}$,
which decreases when $m$ increases, and the variance term, proportional
to $D_{m}/(n\Delta)$. If $g$ belongs to the Besov space $\rond{B}_{2,\infty}^{\alpha}$,
then the bias term $\norm{g_{m}-g_{A}}_{L^{2}}^{2}\propto D_{m}^{-2\alpha}$.
The risk is then minimum for $m_{opt}=(n\Delta)^{1/(1+2\alpha)}$,
and satisfies \[
\rond{R}_{n}(\hat{g}_{m_{opt}})\lesssim(n\Delta)^{-2\alpha/(2\alpha+1)}+\Delta.\]

\subsection{Adaptive estimator}

To bound the risk of the adaptive estimator, we need the additional
assumption:~

\begin{assumption}\label{hypo_nu_sous_exp} 
\begin{enumerate}
\item The Lévy measure $\nu$ is sub-exponential: \[
\exists\lambda,C>0,\quad\forall\vert z\vert>1,\quad\nu(]-z,z[^{c})\leq Ce^{-\lambda\left|z\right|}.\]

\item There exists $\eta$, $\eta>1$, such that $\Delta^{\eta}=O(n^{-1})$. 
\end{enumerate}
\end{assumption}

Let us consider the penalty function $pen(m)=\kappa\xi_{0}^{4}\frac{D_{m}}{n\Delta}$
and choose the adaptive estimator $\hat{g}_{\hat{m}}$ by minimizing
the function

\[
\hat{m}=\min_{m\in\rond{M}_{n}}\gamma_{n}(\hat{g}_{m})+pen(m).\]

We introduce the function $p(m,m')=\frac{pen(m)+pen(m')}{12}.$ For
any $m\in\rond{M}_{n}$,

\begin{align*}
\E{\norm{\hat{g}_{\hat{m}}-g_{A}}_{n}^{2}\units{\Omega_{n}}} & \lesssim\norm{g_{m}-g}_{L^{2}}^{2}+\E{A_{k\Delta}^{2}}+2pen(m)\\
 & +12\E{\sum_{m'\in\rond{M}_{n}}\left[\left(\sup_{t\in\rond{B}_{m,m'}}\nu_{n}^{2}(t)-p(m,m')\right)\units{\Omega_{n}}\right]_{+}}\end{align*}
 where $\rond{B}_{m,m'}=\left\{ t\in S_{m}+S_{m'},\norm{t}_{\pi}\leq1\right\} $.
In order to bound the remaining term, \[
\E{\left[\sup_{t\in\rond{B}_{m,m'}}\nu_{n}^{2}(t)-p(m,m')\right]_{+}},\]
 we use the Berbee's coupling Lemma and a Talagrand's inequality.
Berbee's coupling Lemma is proved by \citet{viennet97}. As the random
variables $(X_{k\Delta})$ are exponentially $\beta-$mixing, it allows
us to deal with independent random variables. 

\begin{berb}

Let $(X_{t})_{t\geq0}$ be a stationary and exponentially $\beta-$mixing
process observed at discrete times $t=0,\Delta,\ldots,n\Delta$. Let
us set $n=2p_{n}q_{n}$ with $q_{n}=8\ln(n)/\Delta$. For any $a\in\left\{ 0,1\right\} $,
$1\leq k\leq p_{n}$, we consider the random variables \[
U_{k,a}=\left(X_{((2(k-1)+a)q_{n}+1)\Delta},\ldots,X_{(2k-1+a)q_{n}\Delta}\right).\]
 There exist random variables $X_{\Delta}^{*},\ldots,X_{n\Delta}^{*}$
such that \[
U_{k,a}^{*}=\left(X_{((2(k-1)+a)q_{n}+1)\Delta}^{*},\ldots,X_{(2k-1+a)q_{n}\Delta}^{*}\right)\]
 satisfy: 
\begin{itemize}
\item For any $a\in\{0,1\}$, the random vectors $U_{1,a}^{*},U_{2,a}^{*},\ldots,U_{p_{n},a}^{*}$
are independent. 
\item For any $(a,k)\in\{0,1\}\times\left\{ 1,\ldots,p_{n}\right\} $, $U_{k,a}^{*}\sim U_{k,a}$. 
\item For any $(a,k)\in\{0,1\}\times\{1,\ldots,p_{n}\}$, $\mathbb{P}\left(U_{k,a}\neq U_{k,a}^{*}\right)\leq\beta(q_{n}\Delta)\leq n^{-8}$. 
\end{itemize}
Let us set $\Omega^{*}=\left\{ \omega,\:\forall(k,a)\in\{0,1\}\times\{1,\ldots,p_{n}\},\: U_{k,a}=U_{k,a}^{*}\right\} $.
Then $\mathbb{P}(\Omega^{*})\leq n\Delta/n^{8}$.

\end{berb}

The following Talagrand's inequality is proved by \citet{birgemassart98}
(corollary 2p.354) and \citet{comtemerlevede2002} (p222-223). 

\begin{talagrand}

Let $(X_{1},\ldots,X_{n})$ be independent identically distributed
random variables and $f_{n}:\rond{B}_{m,m'}\rightarrow S_{m}$ such
that \[
f_{n}(t)=\frac{1}{n}\sum_{k=1}^{n}F_{t}(X_{k})-\E{F_{t}(X_{k})}.\]
 If \[
\sup_{t\in\rond{B}_{m,m'}}\norm{F_{t}}_{\infty}\leq M,\quad\E{\sup_{t\in\rond{B}_{m,m'}}f_{n}^{2}(t)}\leq H^{2},\quad\sup_{t\in\rond{B}_{m,m'}}\var{F_{t}(X_{k})}\leq V\]
 then \[
\E{\sup_{t\in\rond{B}_{m,m'}}f_{n}^{2}(t)-12H^{2}}_{+}\lesssim\frac{V}{n}\exp\left(-k_{1}\frac{nH^{2}}{V}\right)+\frac{M^{2}}{n^{2}}\exp\left(-k_{2}\frac{nH}{M}\right).\]

\end{talagrand}

We then obtain the following oracle inequality:

\begin{thme}\label{thme_est_sigma_xi_adaptatif}

Under assumptions A\ref{hypo_modele}-A\ref{hypo_nu_sous_exp}, there
exists $\kappa_{0}$ such that for any $\kappa\geq\kappa_{0}$, \[
\rond{R}_{n}\left(\hat{g}_{\hat{m}}\right)\lesssim\inf_{m\in\rond{M}_{n}}\left\{ \norm{g_{m}-g_{A}}_{L^{2}}^{2}+pen(m)\right\} +\Delta+\frac{\ln^{3}(n)}{n\Delta}.\]

\end{thme}

The adaptive estimator $\hat{g}_{\hat{m}}$ automatically realises
the best (up to a multiplicative constant) compromise.

\section{Estimation of $\sigma^{2}$. \label{sec:Estimation-sigma}}

We have that \[
T_{k\Delta}=\frac{(X_{(k+1)\Delta}-X_{k\Delta})^{2}}{\Delta}=\sigma^{2}(X_{k\Delta})+\frac{1}{\Delta}J_{k\Delta}^{2}+\textrm{small terms}+\textrm{ centred terms}.\]
 The idea is to keep $T_{k\Delta}$ only when there is no jumps. As
the stochastic term $Z_{k\Delta}$ is of order $\Delta^{1/2}$, we
can only suppress the jumps of amplitude greater than $\Delta^{1/2}$.
Then we consider:

\[
Y_{k\Delta}=\frac{\left(X_{(k+1)\Delta}-X_{k\Delta}\right)^{2}}{\Delta}\units{\Omega_{X,k}}\]
 where $\Omega_{X,k}=\left\{ \omega,\:\left|X_{(k+1)\Delta}-X_{k\Delta}\right|\leq\left(\sigma_{0}+\xi_{0}\right)\ln(n)\Delta^{1/2}+\Delta^{1/2}\right\} $.
We have that \begin{align*}
Y_{k\Delta} & =\sigma^{2}(X_{k\Delta})-\sigma^{2}(X_{k\Delta})\units{\Omega_{X,k}^{c}}+\xi^{2}(X_{k\Delta^{-}})\units{\Omega_{X,k}}+\left(A_{k\Delta}+B_{k\Delta}+E_{k\Delta}\right)\units{\Omega_{X,k}}\\
 & =\sigma^{2}(X_{k\Delta})-\sigma^{2}(X_{k\Delta})\units{\Omega_{X,k}^{c}}+\left(\tilde{A}_{k\Delta}+B_{k\Delta}+\tilde{E}_{k\Delta}\right)\units{\Omega_{X,k}}\end{align*}
 with $\tilde{A}_{k\Delta}=A_{k\Delta}^{(1)}+A_{k\Delta}^{(2)}+A_{k\Delta}^{(3)}$
and $\tilde{E}_{k\Delta}=E_{k\Delta}^{(1)}+E_{k\Delta}^{(2)}+\frac{1}{\Delta}\left(\int_{k\Delta}^{(k+1)\Delta}\xi(X_{s^{-}})dL_{s}\right)^{2}$.
Let us consider $J_{k\Delta}^{(i)}=\int_{k\Delta}^{(k+1)\Delta}\xi(X_{s}^{-})dL_{s}^{(i)}$,
with \begin{gather*}
L_{s}^{(1)}=\int_{\vert z\vert\leq\Delta^{1/2}}z\mu(dz,ds),\quad L_{s}^{(2)}=\int_{\vert z\vert\in]\Delta^{1/2},\Delta^{1/4}]}z\mu(dz,ds),\\
L_{s}^{(3)}=\int_{\vert z\vert>\Delta^{1/4}}z\mu(dz,ds)\end{gather*}
 and denote by $N_{k}=\mu\left(\left](k\Delta,(k+1)\Delta\right],\left[-\Delta^{1/4},\Delta^{1/4}\right]^{c}\right)$
the number of jumps of amplitude greater than $\Delta^{1/4}$ on the
time interval $]k\Delta,(k+1)\Delta]$. We introduce the set \[
\Omega_{N,k}=\left\{ \omega,\: N_{k}=0\quad\textrm{and}\quad\left|\int_{k\Delta}^{(k+1)\Delta}dL_{s}^{(1)}+dL_{s}^{(2)}\right|\leq4\frac{\sigma_{0}+\xi_{0}}{\xi_{1}}\Delta^{1/2}\ln(n)\right\} .\]
The term $B_{k\Delta}\units{\Omega_{X,k}}$ is no longer centred.
Let us set \[
\tilde{B}_{k\Delta}=B_{k\Delta}\units{\Omega_{X,k}\cap\Omega_{N,k}}-\E{B_{k\Delta}\units{\Omega_{X,k}\cap\Omega_{N,k}}\left|\rond{F}_{k\Delta}\right.}\]
 and \begin{align*}
F_{k\Delta} & =\left(\tilde{A}_{k\Delta}+\tilde{E}_{k\Delta}\right)\units{\Omega_{X,k}}-\sigma^{2}(X_{k\Delta})\units{\Omega_{X,k}^{c}}+B_{k\Delta}\units{\Omega_{X,k}\cap\Omega_{N,k}^{c}}\\
 & -\E{B_{k\Delta}\units{\Omega_{X,k}\cap\Omega_{N,k}}\left|\rond{F}_{k\Delta}\right.}.\end{align*}
Then \[
Y_{k\Delta}=\sigma^{2}(X_{k\Delta})+F_{k\Delta}+\tilde{B}_{k\Delta}.\]
 The following assumption is needed.

\begin{assumption}\label{hypo_xi_bornee} 
\begin{enumerate}
\item The function $\xi$ is bounded from below: $\exists\xi_{1}$, $\forall x\in\mathbb{R}$,
$\xi^{2}(x)\geq\xi_{1}^{2}>0$. 
\item There exists $\eta$, $\eta>1$ such that $\Delta^{\eta}=O(n^{-1})$. 
\end{enumerate}
\end{assumption}

The following lemmas are proved later. 

\begin{lem}\label{lem_majoration_omega_X}

$\mathbb{P}\left(\Omega_{X,k}^{c}\right)\lesssim\Delta^{1-\beta/2}+n^{-1}$,
$\mathbb{P}(\Omega_{N,k}^{c})\lesssim\Delta^{1-\beta/2}$ and $\mathbb{P}\left(\Omega_{X,k}\cap\Omega_{N,k}^{c}\right)\lesssim\Delta^{2-\beta/2}+n^{-1}$.

\end{lem}

\begin{lem}\label{lem_majoration_AB-1} 
\begin{itemize}
\item $\E{\tilde{A}_{k\Delta}^{2}\left|\rond{F}_{k\Delta}\right.}\lesssim\Delta$
and $\E{\tilde{A}_{k\Delta}^{4}\left|\rond{F}_{k\Delta}\right.}\lesssim\Delta$. 
\item $\E{\tilde{B}_{k\Delta}\left|\rond{F}_{k\Delta}\right.}=0$, $\E{\tilde{B}_{k\Delta}^{2}\left|\rond{F}_{k\Delta}\right.}\leq\sigma_{0}^{4}/n$
 and $\E{\tilde{B}_{k\Delta}^{4}\left|\rond{F}_{k\Delta}\right.}\lesssim1$. 
\item $\E{\tilde{E}_{k\Delta}^{2}\units{\Omega_{X,k}}\left|\rond{F}_{k\Delta}\right.}\lesssim\Delta^{1-\beta/2}$
and $\E{\tilde{E}_{k\Delta}^{4}\units{\Omega_{X,k}}\left|\rond{F}_{k\Delta}\right.}\lesssim1$. 
\end{itemize}
\end{lem}

\subsection{Estimator for fixed $m$}

We consider the following contrast function and the empirical risk
\[
\tilde{\gamma}_{n}(t)=\frac{1}{n}\sum_{k=1}^{n}\left(t(X_{k\Delta})-Y_{k\Delta}\right)^{2}\units{X_{k\Delta}\in A}\quad\textrm{and}\quad\rond{R}_{n}(t)=\E{\norm{t-\sigma^{2}}_{n}^{2}}.\]
 Let us set $\hat{\sigma}_{m}^{2}=\arg\inf_{t\in S_{m}}\tilde{\gamma}_{n}(t)$.

\begin{thme}\label{thme_sigma_m_fixe}

Under Assumptions A\ref{hypo_modele}, A\ref{hypo_stationnaire.},
A\ref{hypo_espaces} and A\ref{hypo_xi_bornee}, we have that

\[
\rond{R}_{n}(\hat{\sigma}_{m}^{2})\lesssim\left\Vert \sigma_{A}^{2}-\sigma_{m}^{2}\right\Vert _{L^{2}}^{2}+\sigma_{0}^{4}\frac{D_{m}}{n}+\Delta^{1-\beta/2}\ln^{4}(n)\]
 where $\sigma_{A}^{2}(x)=\sigma^{2}(x)\units{x\in A}$.

\end{thme}

The bias term $\left\Vert \sigma_{A}^{2}-\sigma_{m}^{2}\right\Vert _{L^{2}}^{2}$
and the variance term $\sigma_{0}^{4}D_{m}n^{-1}$ are the same as
for a diffusion without jumps. Nevertheless, the remainder term is
$\Delta^{2}$ for a diffusion process (see for instance \citet{comtegenon2007}).
Even for Poisson processes, the remainder term will be here proportional
to $\Delta\ln^{4}(n)$.

If $\sigma^{2}$ belongs to $\rond{B}_{2,\infty}^{\alpha}$, then
$\norm{\sigma_{A}^{2}-\sigma_{m}^{2}}_{L^{2}}^{2}\lesssim D_{m}^{-2\alpha}$
. The best estimator is obtained for $D_{m_{opt}}=n^{-1/(1+2\alpha)}$
and its risk is bounded by $n^{-2\alpha/(2\alpha+1)}+\Delta^{1-\beta/2}$.

\begin{req}\label{req_rate_convergence}

Let us set $\Delta\sim n^{-a}$, with $0<a<1$. We have the following
rates of convergence:

\begin{tabular}{|c|c|c|}
\hline 
$a$  & jumps diffusions  & diffusions\tabularnewline
\hline
\hline 
$0<a\leq\frac{2\alpha}{2(2\alpha+1)}\leq\frac{1}{2}$  & $\Delta^{1/2-\beta/4}$  & $\Delta$\tabularnewline
\hline 
$\frac{2\alpha}{2(2\alpha+1)}\leq a\leq\frac{2\alpha}{(2\alpha+1)(1-\beta/2)}\wedge1$  & $\Delta^{1/2-\beta/4}$  & $n^{-\alpha/(2\alpha+1)}$\tabularnewline
\hline 
$\frac{2\alpha}{(2\alpha+1)(1-\beta/2)}\wedge1\leq a<1$  & $n^{-\alpha/(2\alpha+1)}$  & $n^{-\alpha/(2\alpha+1)}$\tabularnewline
\hline
\end{tabular}

If $\beta=0$, the adaptive estimator will reach the rate of convergence
$n^{-\alpha/(2\alpha+1)}$ for high frequency data ($n\Delta^{\text{2\ensuremath{\alpha}}+1/(2\alpha)}=O(1))$.
This is the minimax rate of convergence for non-parametric estimation
of $\sigma^{2}$ for diffusions processes (see for instance \citet{hofmann99}).
If $\beta$ or $\alpha$ is too big (as soon as $\beta(\alpha+1/2)>1$),
even for high frequency data, the remainder term will be predominant
in the risk.

\end{req}

\subsection{Adaptive estimator}

Let us introduce a penalty function $\widetilde{pen}(m)=\kappa n^{-1}\sigma_{0}^{2}$
and define the adaptive estimator $\hat{\sigma}_{\hat{m}}^{2}$: \[
\hat{m}=\arg\min_{m\in\rond{M}_{n}}\tilde{\gamma}_{n}\left(\hat{\sigma}_{m}^{2}\right)+\widetilde{pen}(m)\]
 where $\rond{M}_{n}=\left\{ m,\: D_{m}\leq\rond{D}_{n}\right\} $.
As for the adaptive estimator of $g=\sigma^{2}+\xi^{2}$, we use the
Berbee's coupling lemma and the Talagrand's inequality to bound the
risk of the estimator $\hat{\sigma}_{\hat{m}}$.

\begin{thme}\label{thme_sigma_adaptatif}

Under Assumptions A\ref{hypo_modele}, A\ref{hypo_stationnaire.},
A\ref{hypo_espaces} and A\ref{hypo_xi_bornee}, there exists $\kappa_{1}$
such that, if $\kappa\geq\kappa_{1}$, we have the following oracle
inequality: \[
\rond{R}_{n}(\hat{\sigma}_{\hat{m}}^{2})\lesssim\min_{m\in\rond{M}_{n}}\left(\norm{\sigma_{A}^{2}-\sigma_{m}^{2}}_{L^{2}}^{2}+\widetilde{pen}(m)\right)+\Delta^{1-\beta/2}\ln^{4}(n)+\frac{1}{n}.\]

\end{thme}

\begin{req}

If Assumptions A\ref{hypo_modele}-A\ref{hypo_xi_bornee} are satisfied,
the risk of the estimator $\hat{\xi}^{2}=\hat{g}_{\hat{m}}-\hat{\sigma}_{\hat{m}}^{2}$satisfies
the following inequality:~ \begin{eqnarray*}
\E{\norm{\hat{\xi}^{2}-\xi_{A}^{2}}} & \lesssim & \min_{m\in\rond{M}_{n}}\left\{ \norm{g_{m}-g_{A}}_{\pi}^{2}+\kappa_{0}\frac{D_{m}}{n\Delta}\right\} +\min_{m\in\rond{M}_{n}}\left\{ \norm{\sigma_{m}^{2}-\sigma_{A}^{2}}+\kappa_{1}\frac{D_{m}}{n}\right\} \\
 & + & \Delta^{1-\beta/2}\ln^{2}(n).\end{eqnarray*}

\end{req}

\section{Simulations \label{sec:Simulations}}

\subsection{Models}

We consider a stochastic process $(X_{t})$ such that \[
dX_{t}=b(X_{t})dt+\sigma(X_{t})dW_{t}+\xi(X_{t^{-}})dL_{t},\quad X_{0}=\eta,\]
 with $L_{t}$ a compound Poisson process:~ \[
L_{t}=\sum_{k=1}^{N_{t}}\zeta_{k}\]
 where $N_{t}$ is a compound Poisson process of intensity 1, and
$(\zeta_{k})$ are centred, independent, and identically distributed
random variables. We denote by $F$ the law of $\zeta$ and we assume
that $\E{\zeta_{k}^{2}}=1$ and that the random variables $(\zeta_{k})$
are independent of $(\eta,\left(W_{t}\right)_{t\geq0},N_{t})$.

\subsubsection{Model 1: Ornstein Uhlenbeck}

\[
dX_{t}=-2X_{t}dt+dW_{t}+dL_{t}\]
 with binomial jumps: $\mathbb{P}(\zeta=1)=\mathbb{P}(\zeta=-1)=0.5$.

\subsubsection{Model 2~}

\[
dX_{t}=-2X_{t}dt+\frac{X_{t^{-}}^{2}+3}{X_{t^{-}}^{2}+1}dW_{t}+dL_{t}\]
 with Laplace jumps: \[
f(dz)=\nu(dz)=0.5e^{-\lambda\vert x\vert}.\]

\subsubsection{Model 3}

\[
dX_{t}=(-2X_{t}+\sin(3X_{t}))dt+\sqrt{2+0.5\sin(\pi X_{t^{-}})}(dW_{t}+dL_{t})\]
 with normal jumps: $\zeta_{k}\sim\rond{N}(0,1)$.

\subsubsection{Model 4:~}

In this model, the Lévy process is not a compound Poisson process.
We set \[
n(z)=\sum_{k=1}^{\infty}2^{k+1}(\delta_{1/2^{k}}+\delta_{-1/2^{k}}),\quad b(x)=-2x\quad\textrm{and}\quad\sigma(x)=\xi(x)=1.\]
 The Blumenthal-Getoor index of this process is such that $\beta>1$.

\subsection{Method}

We use the vectorial subspaces generated by the spline functions:
\begin{gather*}
S_{m,r}=\textrm{Vect}\left(\varphi_{r,k,m},\; k\in\mathbb{Z}\right),\quad\textrm{with}\quad\varphi_{r,k,m}=2^{m/2}g_{r}(2^{m}x-k)\units{x\in A}\\
\textrm{and}\quad g_{r}=\units{x\in A}*\ldots*\units{x\in A}\end{gather*}
 Those subspaces form a multi-resolution analysis of $L^{2}(A)$.
We use the same simulation method as in \citet{Rubenthaler_HDR}. 

To construct the adaptive estimator, we compute $\hat{f}_{m,r}$ for
$D_{m}\leq\sqrt{n\Delta}$, $0\leq r\leq4$ and $m\leq7$ (for $m=7$,
we already have $D_{m}=128$. If $m$ was bigger, there will be a
memory problem). Then we minimize $\gamma_{n}(\hat{f}_{m,r})+pen(m,r)$
with respect to $m$, then $r$. There is three constants in the penalty
function $pen(m,r)$. The constants $\sigma_{0}^{4}$ and $\xi_{0}^{4}$
are unknown, but they can be replaced by rough estimators, as only
an upper bound for $\sigma_{0}^{4}$ and $\xi_{0}^{4}$ is needed.
In our simulations, we took the true value of $\sigma_{0}^{4}$ and
$\xi_{0}^{4}$. The constants $\kappa_{0}$ and $\kappa_{1}$ are
chosen by numerical calibration (see \citet{comteroz2002,comteroz2004}
for a complete discussion). Another way of dealing with the constants
of the penalty would be the slope method developed by \citet{arlotmassart2009},
however, this method is a bit slow.

\begin{figure}[p]
 \caption{Model 1}

\label{Flo:figure1} \[
dX_{t}=-2X_{t}dt+dW_{t}+dL_{t},\quad\textrm{binomial jumps}\]

\begin{tabular}{cc}
Estimation of $\sigma^{2}$  & Estimation of $\sigma^{2}+\xi^{2}$\tabularnewline
\includegraphics[scale=0.45]{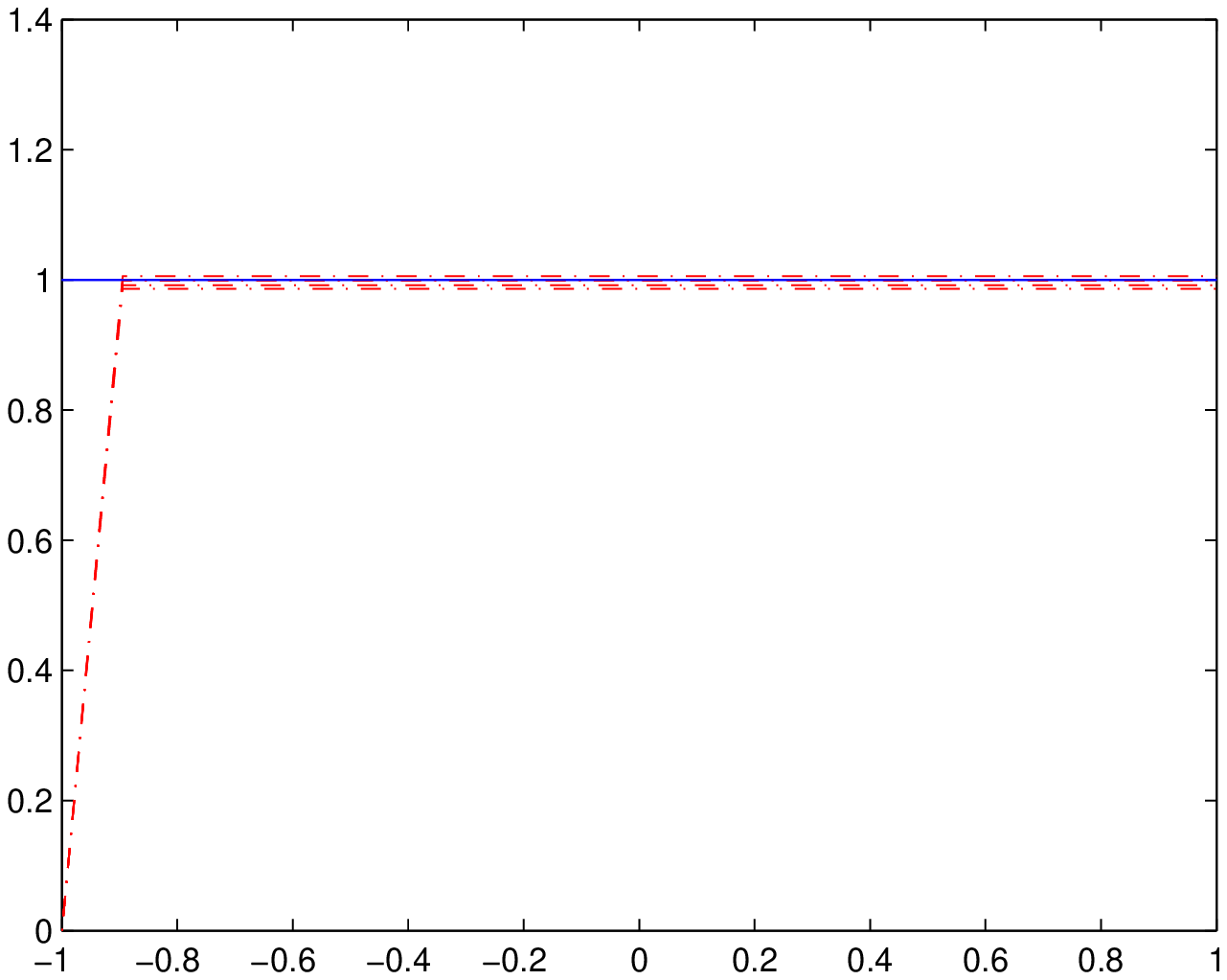}  & \includegraphics[scale=0.45]{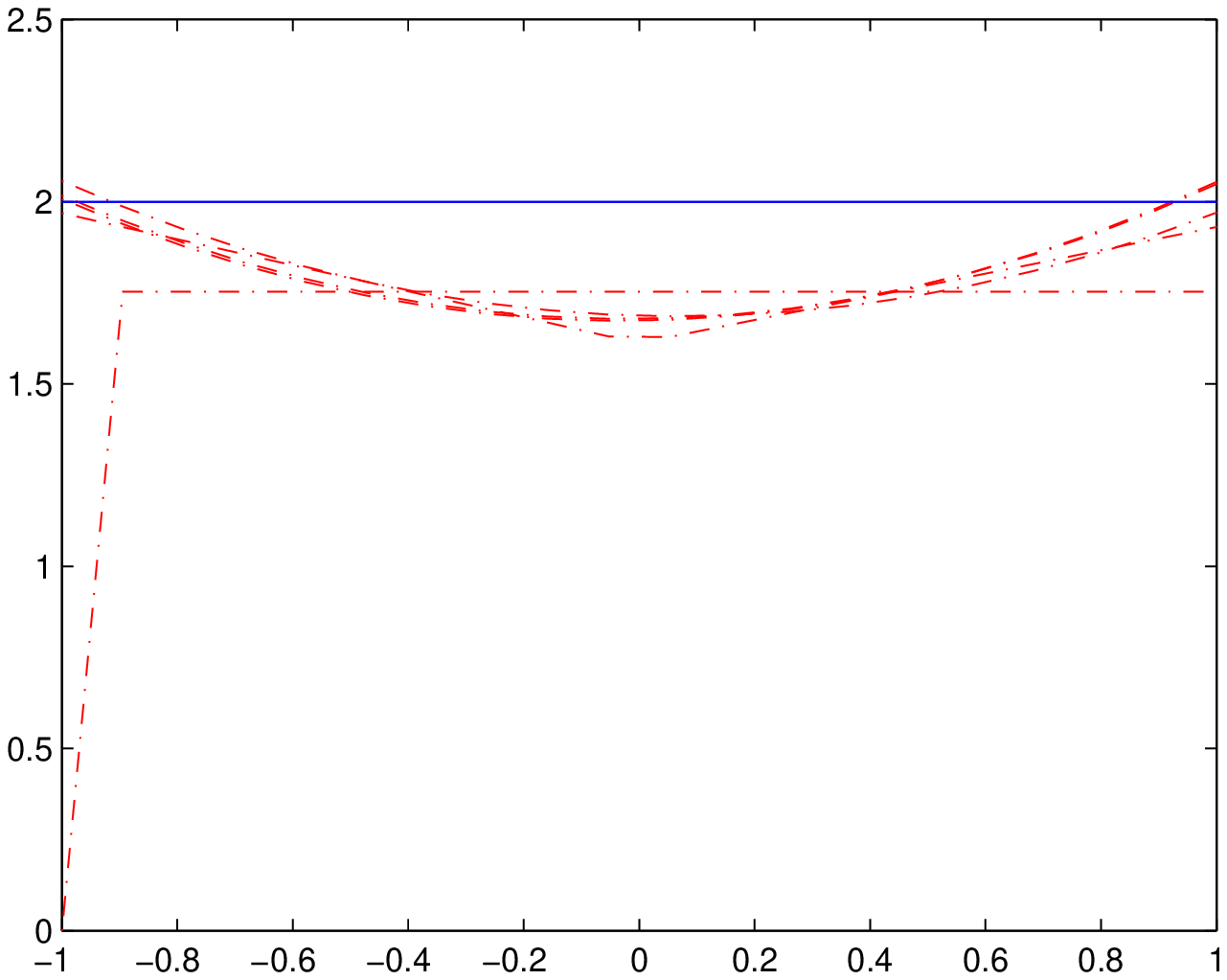}\tabularnewline
$-$ :~true function  & $-$ :~true function\tabularnewline
$-.$ : estimator  & $-.$ : estimator\tabularnewline
$n=10^{5}$, $\Delta=10^{-3}$  & $n=10^{5}$, $\Delta=0.1$ \tabularnewline
\end{tabular}
\end{figure}

\begin{figure}[p]
 \caption{Model 2}

\label{Flo:figure2} \[
dX_{t}=-2X_{t}dt+\frac{X_{t^{-}}^{2}+3}{X_{t^{-}}^{2}+1}dW_{t}+dL_{t},\quad\textrm{Laplace jumps}\]

\begin{tabular}{cc}
Estimation of $\sigma^{2}$  & Estimation of $\sigma^{2}+\xi^{2}$\tabularnewline
\includegraphics[scale=0.45]{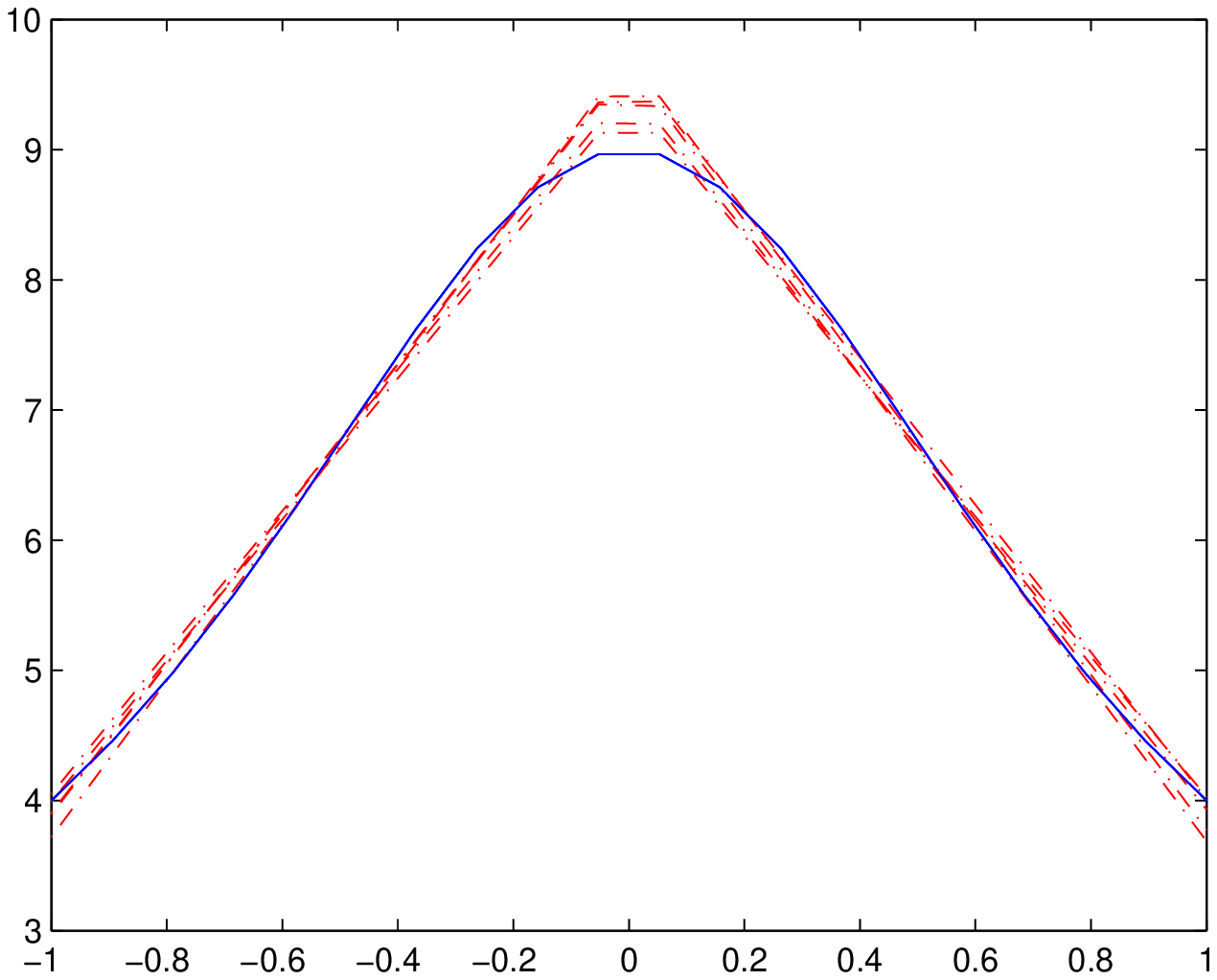}  & \includegraphics[scale=0.45]{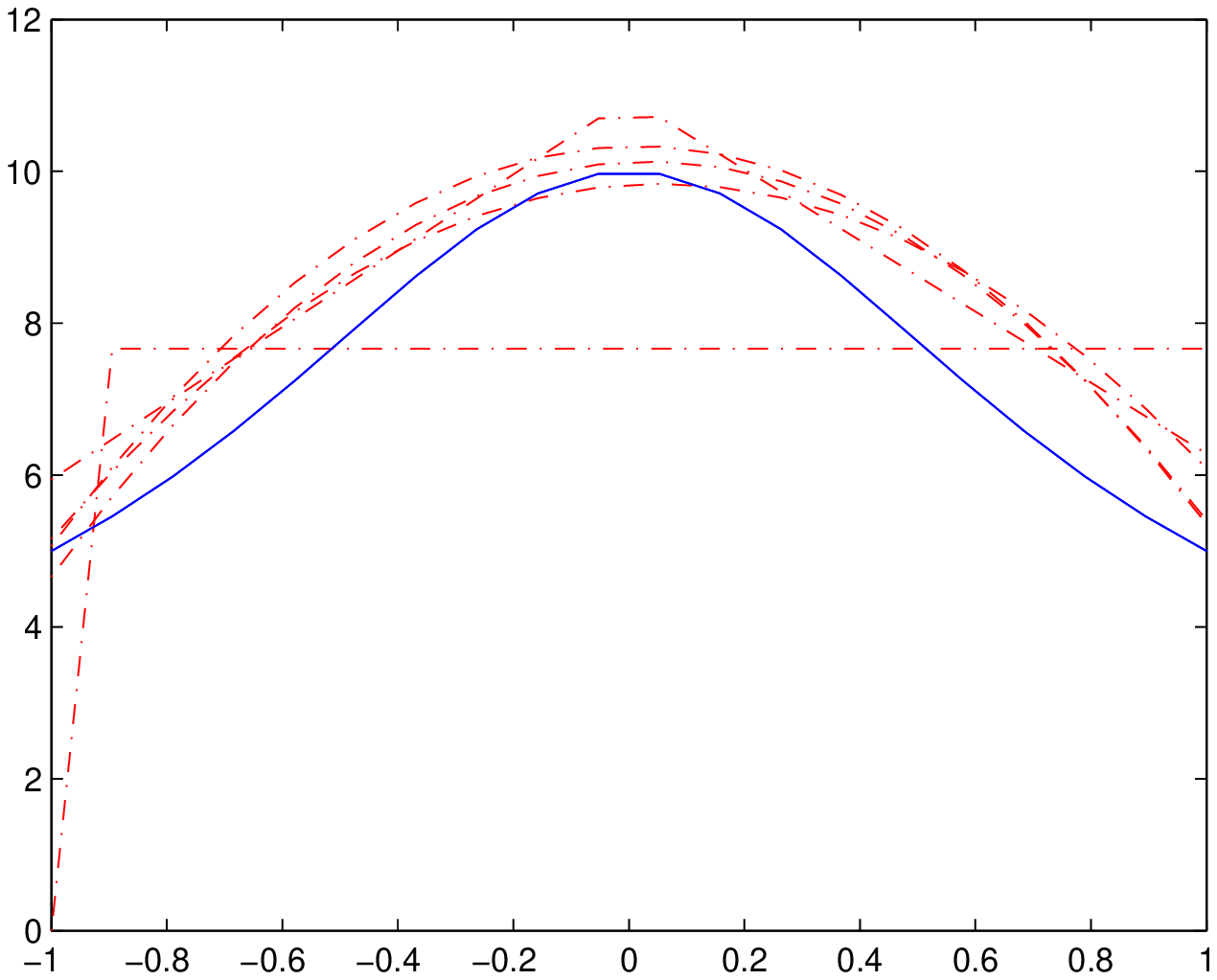}\tabularnewline
$-$ :~true function  & $-$:~true function\tabularnewline
$-.$ : estimator  & $-.$ : estimator\tabularnewline
$n=10^{5}$, $\Delta=10^{-3}$  & $n=10^{5}$, $\Delta=10^{-2}$ \tabularnewline
\end{tabular}
\end{figure}

\begin{figure}[p]
 \caption{Model 3}

\label{Flo:figure3}

\[
dX_{t}=(-2X_{t}+\sin(3X_{t}))dt+\sqrt{2+0.5\sin(\pi X_{t^{-}})}(dW_{t}+dL_{t}),\quad\textrm{Normal jumps}\]

\begin{tabular}{cc}
Estimation of $\sigma^{2}$  & Estimation of $\sigma^{2}+\xi^{2}$\tabularnewline
\includegraphics[scale=0.45]{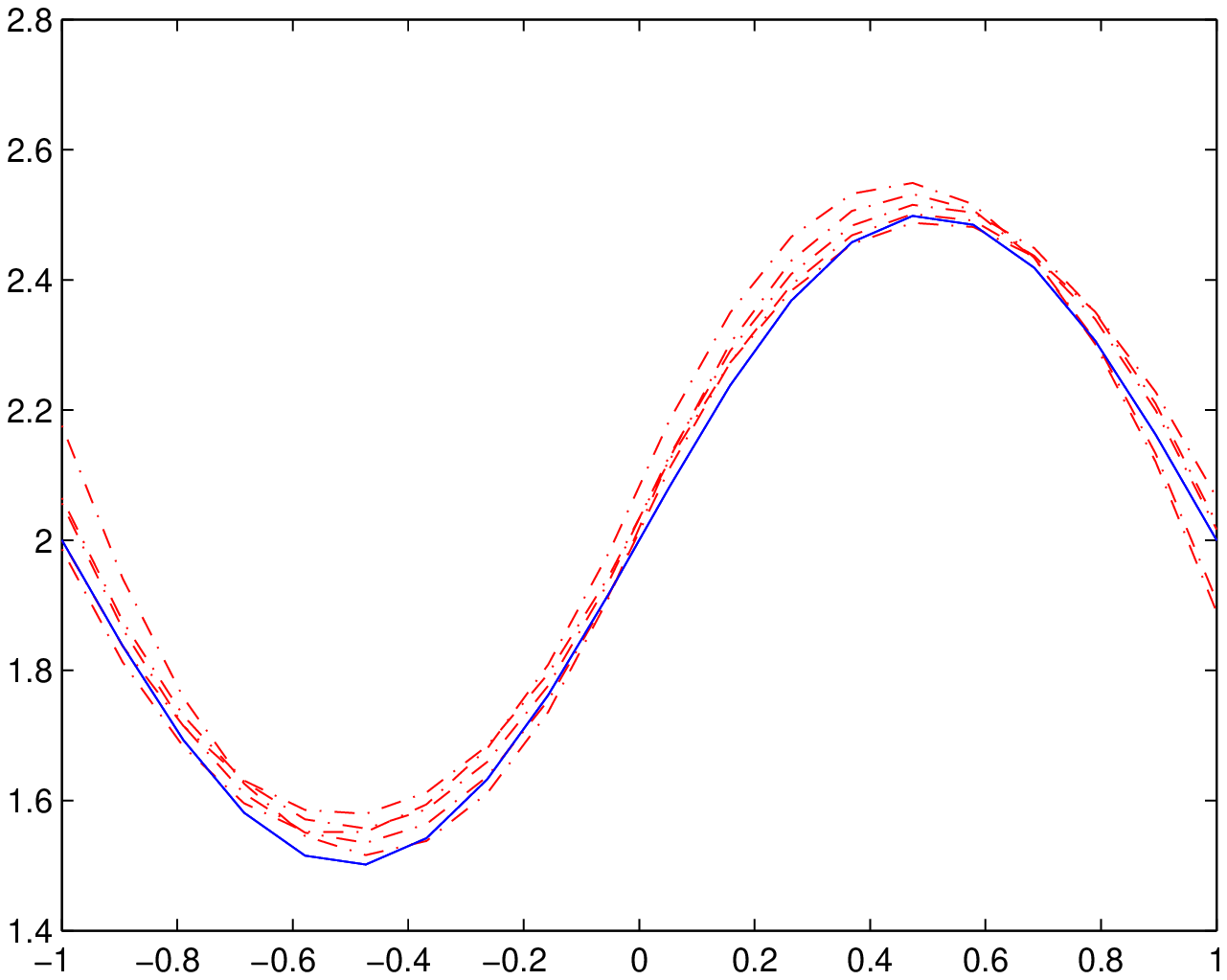}  & \includegraphics[scale=0.45]{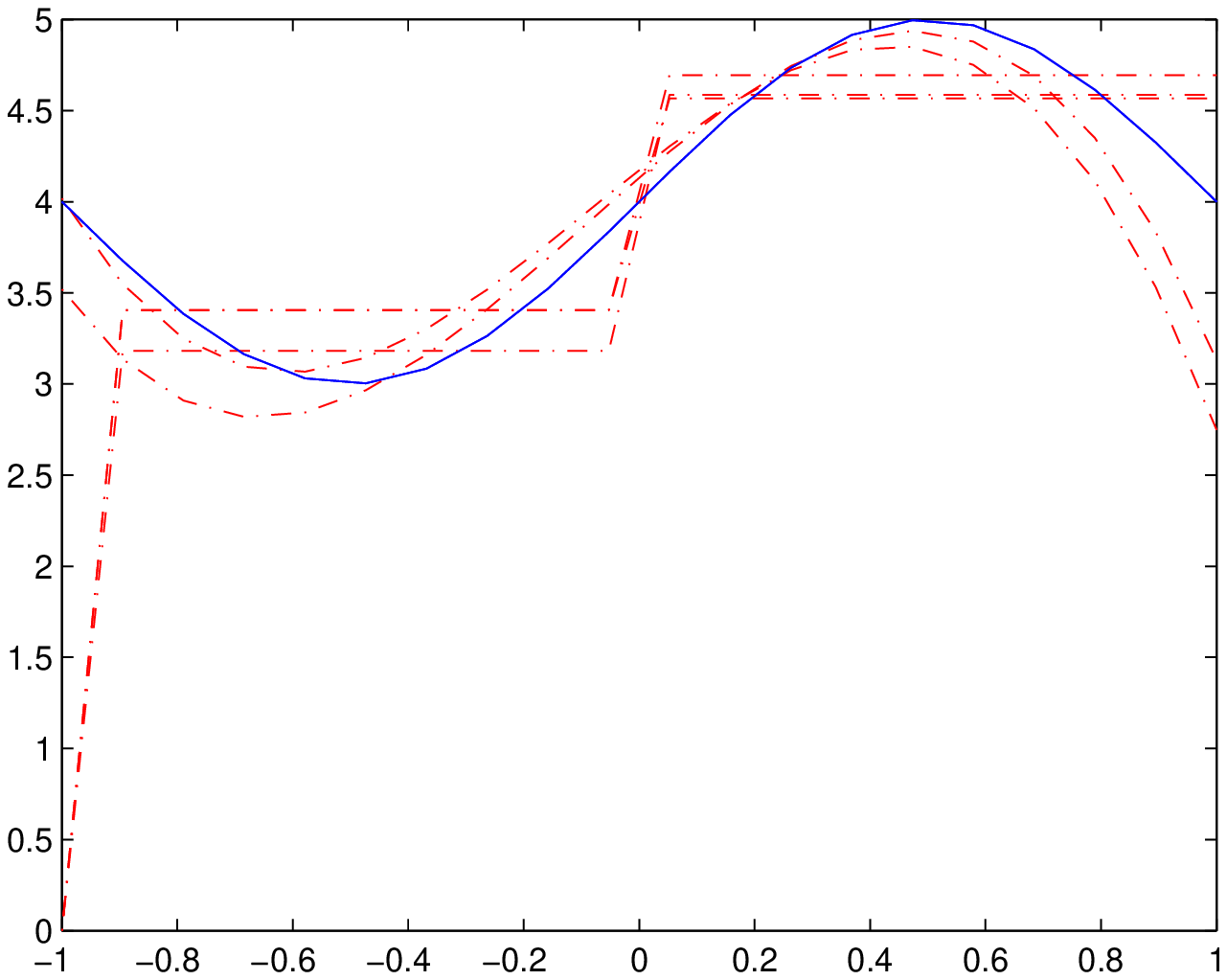}\tabularnewline
$-$ :~true function  & $-$ :~true function\tabularnewline
$-.$ : estimator  & $-.$ : estimator\tabularnewline
$n=10^{5}$, $\Delta=10^{-3}$  & $n=10^{5}$, $\Delta=10^{-2}$ \tabularnewline
\end{tabular}
\end{figure}

To obtain Figures \ref{Flo:figure1}-\ref{Flo:figure4}, for each
model, we realise 5 simulations and draw the 5 corresponding estimators.
To construct Tables \eqref{Flo:table1}-\eqref{Flo:table4}, for each
couplet $(n,\Delta)$ and each model, we make 50 simulations, and
for each simulation, we compute the adaptive estimator $\hat{g}_{\hat{m},\hat{r}}$
or $\hat{\sigma}_{\hat{m},\hat{r}}$ , the selected dimension $\left(\hat{m},\hat{r}\right)$
and the empirical error \begin{gather*}
err=\frac{1}{n}\sum_{k=1}^{n}\left(\hat{g}_{\hat{m},\hat{r}}(X_{k\Delta})-g(X_{k\Delta})\right)^{2}\units{X_{k\Delta}\in A}\\
\textrm{or}\quad err=\frac{1}{n}\sum_{k=1}^{n}\left(\hat{\sigma}_{\hat{m},\hat{r}}^{2}(X_{k\Delta})-\sigma^{2}(X_{k\Delta})\right)^{2}\units{X_{k\Delta}\in A}.\end{gather*}
 We also compute the empirical error for each $\hat{g}_{m,r}$ (or
$\hat{\sigma}_{m,r}^{2}$). Then we deduce the dimension $(m_{min},r_{min})$
that minimizes the empirical error (denoted by $err_{min}$). In the
tables, we write the following informations:~ 
\begin{itemize}
\item mean of the empirical errors of $\hat{g}_{\hat{m},\hat{r}}$ and $\hat{\sigma}_{\hat{m},\hat{r}}^{2}$,
$risk$ 
\item oracle $or=mean(err/err_{min})$. 
\item $m_{est}$ and $r_{est}$, means of $\hat{m}$ and $\hat{r}$. 
\item $t_{e}$ the mean of the estimation time for one simulation. 
\end{itemize}

\subsection{Results}

\begin{figure}[p]
 \caption{Model 4}

\label{Flo:figure4} \[
n(z)=\sum_{k=1}^{\infty}2^{k+1}(\delta_{1/2^{k}}+\delta_{-1/2^{k}}),\quad b(x)=-2x\quad\textrm{and}\quad\sigma(x)=\xi(x)=1.\]

\begin{tabular}{cc}
Estimation of $\sigma^{2}$  & Estimation of $\sigma^{2}+\xi^{2}$\tabularnewline
\includegraphics[scale=0.45]{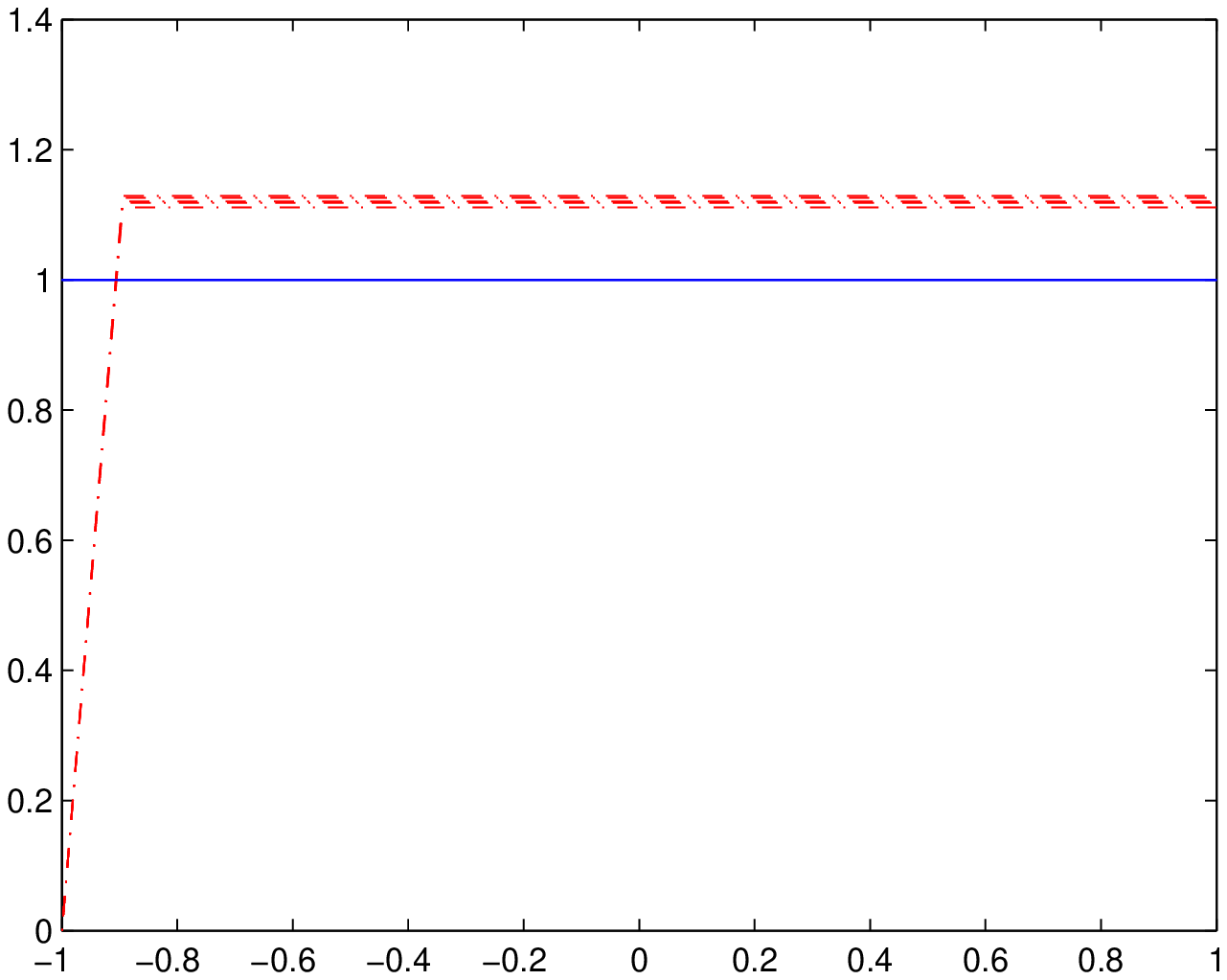}  & \includegraphics[scale=0.45]{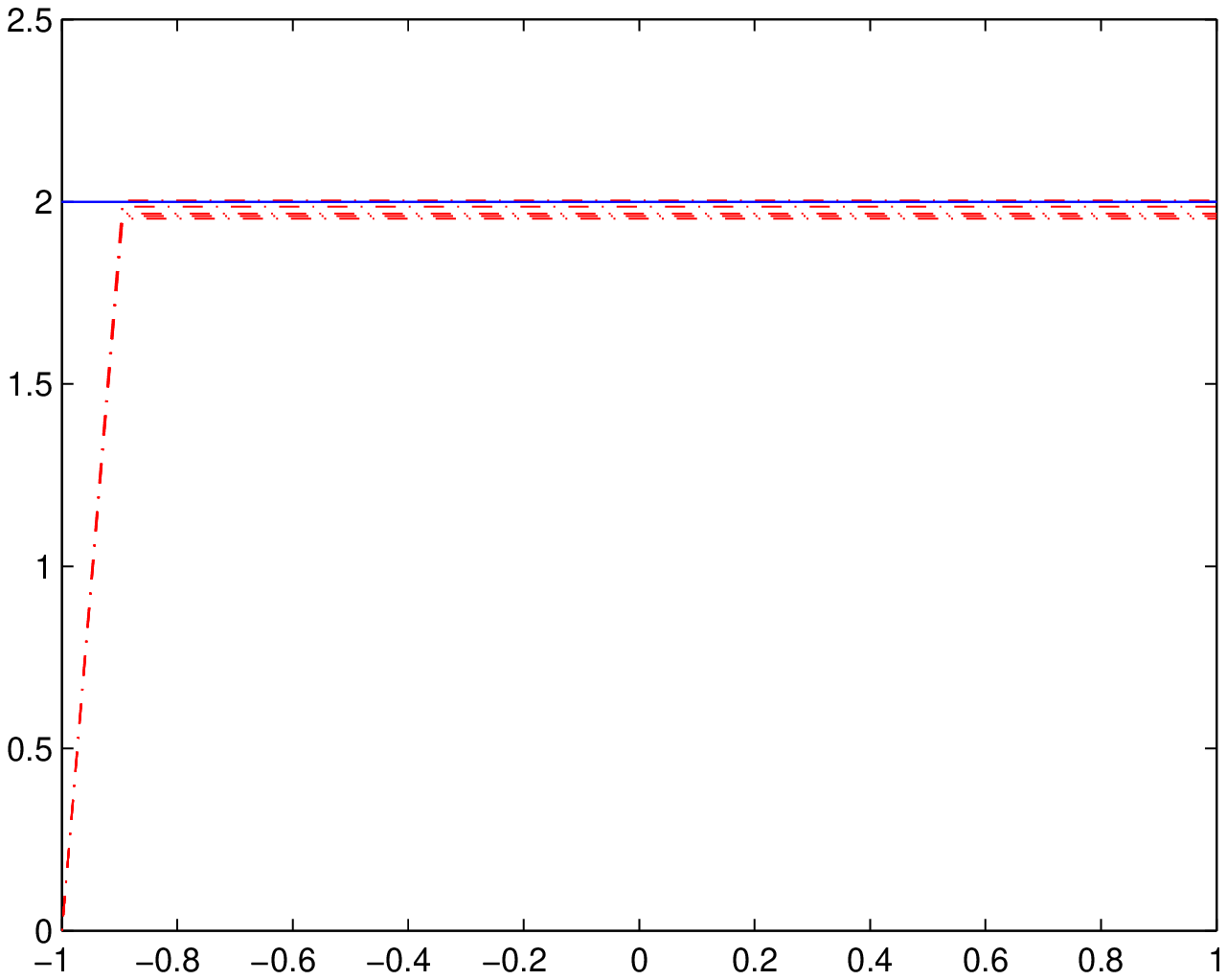}\tabularnewline
$-$ :~true function  & $-$ :~true function\tabularnewline
\_. : estimator  & \_. : estimator\tabularnewline
$n=10^{5}$, $\Delta=10^{-4}$  & $n=10^{5}$, $\Delta=10^{-2}$ \tabularnewline
\end{tabular}
\end{figure}

For Models 1-3, for $\Delta$ small enough ($\Delta=10^{-2}$ or $10^{-3}$
for Model 1, $\Delta=10^{-3}$ for Models 2 and 3), the risk of the
adaptive estimator $\hat{\sigma}_{\hat{m}}^{2}$ is inversely proportional
to $n$, that is proportional to the variance term. In Table \ref{Flo:table4},
we can see that the risk mostly depends on $\Delta$: the remainder
term is predominant. As the Blumethal-Getoor index $\beta>1$, this
is consistent with Remark \eqref{req_rate_convergence}. We can see
in Figure \ref{Flo:figure4} that $\sigma^{2}$ is overestimated:
this is because the small jumps can not be cut. This bias decreases
with $\Delta$.

The function $g=\sigma^{2}+\xi^{2}$ is more difficult to estimate.
Indeed, the variance term is bigger (it is proportional to $1/n\Delta$
and not $1/n$). For $n\Delta$ not big enough ($n\Delta=1$ or 10),
the results can be quite bad. When $\Delta$ is fixed (and small enough
so that the remainder term is not preponderant), the risk decreases
when $n$ increases. 

\begin{table}[p]
 \label{Flo:table1}\caption{Model 1}

\[
dX_{t}=-2X_{t}dt+dW_{t}+dL_{t},\quad\textrm{binomial jumps}\]

\begin{tabular}{|c|c||c|c|c|c|c||c|c|c|c|c|}
\multicolumn{1}{c}{} & \multicolumn{1}{c}{} & \multicolumn{5}{c}{Estimation of $\sigma^{2}+\xi^{2}$} & \multicolumn{5}{c}{Estimation of $\sigma^{2}$}\tabularnewline
\hline 
$\Delta$  & $n$  & risk  & oracle  & $m_{est}$  & $r_{est}$  & $t_{e}$  & risk  & oracle  & $m_{est}$  & $r_{est}$  & $t_{e}$\tabularnewline
\hline
\hline 
$10^{-1}$  & $10^{3}$  & 0.075  & 1.00  & 0.00  & 0.92  & 0.78  & 0.93  & 1.56  & 1.92  & 0.92  & 0.78\tabularnewline
\hline 
$10^{-1}$  & $10^{4}$  & 0.061  & 1.03  & 0.02  & 1.30  & 3.61  & 0.63  & 1.12  & 3.44  & 1.30  & 3.61\tabularnewline
\hline 
$10^{-1}$  & $10^{5}$  & 0.066  & 1.14  & 0.16  & 1.46  & 36  & 0.59  & 1.03  & 3.46  & 1.46  & 36\tabularnewline
\hline 
$10^{-2}$  & $10^{3}$  & 0.15  & 1.00  & 0.00  & 0.00  & 0.22  & 0.0026  & 1.71  & 0.02  & 0.00  & 0.22\tabularnewline
\hline 
$10^{-2}$  & $10^{4}$  & 0.015  & 1.00  & 0.00  & 0.00  & 3.60  & 0.0004  & 3.27  & 0.02  & 0.00  & 3.60\tabularnewline
\hline 
$10^{-2}$  & $10^{5}$  & 0.0021  & 1.00  & 0.00  & 0.52  & 36  & 0.00048  & 4.70  & 0.12  & 0.52  & 36\tabularnewline
\hline 
$10^{-3}$  & $10^{3}$  & 4.18  & 1.21  & 0.02  & 0.00  & 0.13  & 0.0020  & 1.00  & 0.00  & 0.00  & 0.13\tabularnewline
\hline 
$10^{-3}$  & $10^{4}$  & 0.12  & 1.00  & 0.00  & 0.00  & 0.58  & 0.0002  & 1.00  & 0.00  & 0.00  & 0.58\tabularnewline
\hline 
$10^{-3}$  & $10^{5}$  & 0.013  & 1.00  & 0.00  & 0.02  & 36  & 0.000022  & 1.57  & 0.00  & 0.02  & 36\tabularnewline
\hline
\end{tabular}
\end{table}

\begin{table}[p]
 \label{Flo:table2}\caption{Model 2}

\[
dX_{t}=-2X_{t}dt+\frac{X_{t^{-}}^{2}+3}{X_{t^{-}}^{2}+1}dW_{t}+dL_{t},\quad\textrm{Laplace jumps}\]

\begin{tabular}{|c|c||c|c|c|c|c||c|c|c|c|c|}
\multicolumn{1}{c}{} & \multicolumn{1}{c}{} & \multicolumn{5}{c}{Estimation of $\sigma^{2}+\xi^{2}$} & \multicolumn{5}{c}{Estimation of $\sigma^{2}$}\tabularnewline
\hline 
$\Delta$  & $n$  & risk  & oracle  & $m_{est}$  & $r_{est}$  & $t_{e}$  & risk  & oracle  & $m_{est}$  & $r_{est}$  & $t_{e}$\tabularnewline
\hline
\hline 
$10^{-1}$  & $10^{3}$  & 3.53  & 2.46  & 0.00  & 0.02  & 0.78  & 3.02  & 3.04  & 0.00  & 0.02  & 0.77\tabularnewline
\hline 
$10^{-1}$  & $10^{4}$  & 3.07  & 2.05  & 0.00  & 0.42  & 2.54  & 2.23  & 1.68  & 0.14  & 0.42  & 2.54\tabularnewline
\hline 
$10^{-1}$  & $10^{5}$  & 1.51  & 1.01  & 0.52  & 1.58  & 20.4  & 1.26  & 1.01  & 0.46  & 1.58  & 20.4\tabularnewline
\hline 
$10^{-2}$  & $10^{3}$  & 152  & 5.00  & 0.20  & 0.08  & 0.23  & 2.81  & 10.5  & 0.02  & 0.08  & 0.23\tabularnewline
\hline 
$10^{-2}$  & $10^{4}$  & 3.37  & 5.69  & 0.00  & 1.28  & 2.56  & 0.28  & 1.25  & 0.72  & 1.28  & 2.53\tabularnewline
\hline 
$10^{-2}$  & $10^{5}$  & 1.36  & 1.39  & 0.54  & 1.08  & 20.3  & 0.22  & 1.05  & 1.00  & 1.08  & 20.3\tabularnewline
\hline 
$10^{-3}$  & $10^{3}$  & 1600  & 2.87  & 0.02  & 0.20  & 0.14  & 2.34  & 10.2  & 0.16  & 0.20  & 0.14\tabularnewline
\hline 
$10^{-3}$  & $10^{4}$  & 85  & 3.60  & 0.10  & 1.16  & 0.56  & 0.087  & 1.47  & 0.84  & 1.16  & 0.56\tabularnewline
\hline 
$10^{-3}$  & $10^{5}$  & 4.90  & 6.58  & 0.00  & 1.00  & 20.3  & 0.023  & 3.23  & 1.00  & 1.00  & 20.3\tabularnewline
\hline
\end{tabular}
\end{table}

\begin{table}[p]
 \label{Flo:table3}\caption{Model 3}

\[
dX_{t}=(-2X_{t}+\sin(3X_{t}))dt+\sqrt{2+0.5\sin(\pi X_{t^{-}})}(dW_{t}+dL_{t}),\quad\textrm{Normal jumps}\]

\begin{tabular}{|c|c||c|c|c|c|c||c|c|c|c|c|}
\multicolumn{1}{c}{} & \multicolumn{1}{c}{} & \multicolumn{5}{c}{Estimation of $\sigma^{2}+\xi^{2}$} & \multicolumn{5}{c}{Estimation of $\sigma^{2}$}\tabularnewline
\hline 
$\Delta$  & $n$  & risk  & oracle  & $m_{est}$  & $r_{est}$  & $t_{e}$  & risk  & oracle  & $m_{est}$  & $r_{est}$  & $t_{e}$\tabularnewline
\hline
\hline 
$10^{-1}$  & $10^{3}$  & 1.00  & 1.86  & 0.02  & 1.10  & 0.81  & 15.4  & 7.15  & 4.52  & 1.10  & 0.80\tabularnewline
\hline 
$10^{-1}$  & $10^{4}$  & 0.56  & 1.27  & 0.48  & 1.20  & 3.43  & 3.43  & 1.71  & 5.72  & 1.20  & 3.39\tabularnewline
\hline 
$10^{-1}$  & $10^{5}$  & 0.43  & 1.03  & 0.90  & 1.00  & 31.3  & 2.09  & 1.08  & 6.24  & 1.00  & 31.2\tabularnewline
\hline 
$10^{-2}$  & $10^{3}$  & 24.4  & 28.5  & 0.32  & 0.62  & 0.24  & 2.49  & 8.59  & 1.66  & 0.62  & 0.24\tabularnewline
\hline 
$10^{-2}$  & $10^{4}$  & 0.78  & 2.57  & 0.12  & 1.30  & 3.41  & 1.46  & 4.34  & 4.88  & 1.30  & 3.38\tabularnewline
\hline 
$10^{-2}$  & $10^{5}$  & 0.12  & 3.30  & 0.82  & 1.10  & 31.1  & 0.75  & 1.54  & 6.98  & 1.10  & 31.0\tabularnewline
\hline 
$10^{-3}$  & $10^{3}$  & 82.3  & 3.57  & 0.08  & 0.08  & 0.14  & 0.090  & 5.43  & 0.12  & 0.08  & 0.14\tabularnewline
\hline 
$10^{-3}$  & $10^{4}$  & 13.1  & 2.51  & 0.18  & 1.14  & 0.60  & 0.019  & 4.61  & 0.80  & 1.14  & 0.61\tabularnewline
\hline 
$10^{-3}$  & $10^{5}$  & 0.98  & 2.63  & 0.26  & 2.18  & 31  & 0.0026  & 1.16  & 0.82  & 2.18  & 30.8\tabularnewline
\hline
\end{tabular}
\end{table}

\begin{table}[p]
 \label{Flo:table4}\caption{Model 4}

\[
n(z)=\sum_{k=1}^{\infty}2^{k+1}(\delta_{1/2^{k}}+\delta_{-1/2^{k}}),\quad b(x)=-2x\quad\textrm{and}\quad\sigma(x)=\xi(x)=1.\]

\begin{tabular}{|c|c||c|c|c|c|c||c|c|c|c|c|}
\multicolumn{1}{c}{} & \multicolumn{1}{c}{} & \multicolumn{5}{c}{Estimation of $\sigma^{2}+\xi^{2}$} & \multicolumn{5}{c}{Estimation of $\sigma^{2}$}\tabularnewline
\hline 
$\Delta$  & $n$  & risk  & oracle  & $m_{est}$  & $r_{est}$  & $t_{e}$  & risk  & oracle  & $m_{est}$  & $r_{est}$  & $t_{e}$\tabularnewline
\hline
\hline 
$10^{-1}$  & $10^{3}$  & 0.074  & 1.00  & 0.00  & 0.14  & 0.86  & 0.56  & 1.02  & 0.06  & 0.14  & 0.89\tabularnewline
\hline 
$10^{-1}$  & $10^{4}$  & 0.075  & 1.01  & 0.02  & 1.28  & 4.04  & 0.54  & 1.02  & 0.30  & 1.28  & 4.01\tabularnewline
\hline 
$10^{-1}$  & $10^{5}$  & 0.080  & 1.02  & 0.12  & 1.98  & 37.4  & 0.55  & 1.02  & 0.12  & 1.98  & 37.3\tabularnewline
\hline 
$10^{-2}$  & $10^{3}$  & 0.039  & 1.00  & 0.00  & 0.42  & 0.25  & 0.96  & 1.19  & 0.70  & 0.42  & 0.25\tabularnewline
\hline 
$10^{-2}$  & $10^{4}$  & 0.0040  & 1.00  & 0.00  & 0.62  & 4.57  & 0.86  & 1.01  & 0.72  & 0.62  & 4.58\tabularnewline
\hline 
$10^{-2}$  & $10^{5}$  & 0.0012  & 1.00  & 0.00  & 0.58  & 38.3  & 0.91  & 1.00  & 1.24  & 0.58  & 38.1\tabularnewline
\hline 
$10^{-3}$  & $10^{3}$  & 1.22  & 1258  & 0.04  & 1.02  & 0.14  & 0.071  & 1.07  & 0.04  & 0.02  & 0.15\tabularnewline
\hline 
$10^{-3}$  & $10^{4}$  & 0.012  & 1.00  & 0.00  & 0.10  & 0.87  & 0.094  & 1.01  & 0.02  & 0.10  & 0.87\tabularnewline
\hline 
$10^{-3}$  & $10^{5}$  & 0.0015  & 1.00  & 0.00  & 0.36  & 38.5  & 0.15  & 1.00  & 0.22  & 0.36  & 38.5\tabularnewline
\hline 
$10^{-4}$  & $10^{4}$  & 0.24  & 1.00  & 0.00  & 0.02  & 0.39  & 0.013  & 1.04  & 0.02  & 0.02  & 0.39\tabularnewline
\hline 
$10^{-4}$  & $10^{5}$  & 0.021  & 1.00  & 0.00  & 0.04  & 6.31  & 0.014  & 1.00  & 0.00  & 0.04  & 6.26\tabularnewline
\hline
\end{tabular}
\end{table}

\section{Proofs \label{sec:Proofs}}

\subsection{Proof of Theorem \ref{thme_risque_est_sigma_xi_m_fixe}}

By Lemma \ref{lem_majoration_AB}, $\E{A_{k\Delta}^{2}}\lesssim\Delta$.
It remains to bound $\E{\sup_{t\in\rond{B}_{m}}\nu_{n}^{2}(t)}$.
Let $\left(\varphi_{\lambda}\right)_{1\leq\lambda\leq D_{m}}$ be
an orthonormal (for the $\norm{.}_{\pi}$ norm) basis of $S_{m}$.
Any function $t\in\rond{B}_{m}$ can be written $t=\sum_{\lambda=1}^{D_{m}}a_{\lambda}\varphi_{\lambda}$
with $\sum_{\lambda=1}^{D_{m}}a_{\lambda}^{2}\leq1$. By Cauchy-Schwarz,
\[
\sup_{t\in\rond{B}_{m}}\nu_{n}^{2}(t)=\sup_{\sum_{\lambda=1}^{D_{m}}a_{\lambda}^{2}\leq1}\left(\sum_{\lambda=1}^{D_{m}}a_{\lambda}^{2}\right)\left(\sum_{\lambda=1}^{D_{m}}\nu_{n}^{2}\left(\varphi_{\lambda}\right)\right)\leq\sum_{\lambda=1}^{D_{m}}\nu_{n}^{2}(\varphi_{\lambda}).\]
 According to Lemma \ref{lem_majoration_AB}, $\E{B_{k\Delta}+E_{k\Delta}\left|\rond{F}_{k\Delta}\right.}=0$,
and \begin{align*}
\E{\nu_{n}^{2}(\varphi_{\lambda})} & =\EE{\left(\frac{1}{n}\sum_{k=1}^{n}\left(B_{k\Delta}+E_{k\Delta}\right)\varphi_{\lambda}(X_{k\Delta})\right)^{2}}\\
 & \leq\frac{2}{n^{2}}\sum_{k=1}^{n}\EE{\varphi_{\lambda}^{2}(X_{k\Delta})\E{B_{k\Delta}^{2}+E_{k\Delta}^{2}\left|\rond{F}_{k\Delta}\right.}}\\
 & \lesssim\frac{\xi_{0}^{4}}{n\Delta}+\frac{\sigma_{0}^{4}+\sigma_{0}^{2}\xi_{0}^{2}}{n}.\end{align*}
 Then \begin{equation}
\E{\sup_{t\in\rond{B}_{m}}\nu_{n}^{2}(t)}\lesssim\left(\frac{\xi_{0}^{4}}{\Delta}+\sigma_{0}^{4}+\sigma_{0}^{2}\xi_{0}^{2}\right)\frac{D_{m}}{n}\label{eq:majoration_nu_n}\end{equation}

It remains to bound the risk on $\Omega_{n}^{c}$. ~By Lemma \ref{majoration_omegaC},
$\mathbb{P}\left(\Omega_{n}^{c}\right)\leq1/n^{8}$. The function
$\hat{g}_{m}$ is the orthogonal projection (for the $\norm{.}_{n}$
norm) of $(T_{\Delta},\ldots,T_{n\Delta})$ on the vectorial subspace
$\left\{ (t(X_{\Delta}),\ldots,t(X_{n\Delta}))\:,t\in S_{m}\right\} $.
Let us denote by $\Pi_{m}$ the orthogonal projection on this subspace.
As $T_{k\Delta}=g(X_{k\Delta})+A_{k\Delta}+B_{k\Delta}+E_{k\Delta}$,
we obtain: \begin{align*}
\norm{\hat{g}_{m}-g_{A}}_{n}^{2} & =\norm{\Pi_{m}T-g_{A}}_{n}^{2}=\norm{\Pi_{m}g_{A}-g_{A}}_{n}^{2}+\norm{\Pi_{m}A+\Pi_{m}B+\Pi_{m}E}_{n}^{2}\\
 & \leq\norm{g_{A}}_{n}^{2}+\norm{A+B+E}_{n}^{2}\end{align*}
 By stationarity and Cauchy-Schwarz:~ \begin{align*}
\E{\norm{\hat{g}_{m}-g_{A}}_{n}^{2}\units{\Omega_{n}^{c}}} & \lesssim\E{\norm{g_{A}}_{n}^{2}\units{\Omega_{n}^{c}}}+\EE{\left(\frac{1}{n}\sum_{k=1}^{n}A_{k\Delta}^{2}+B_{k\Delta}^{2}+E_{k\Delta}^{2}\right)\units{\Omega_{n}^{c}}}\\
 & \lesssim\left[\left(\EE{\norm{g_{A}}_{n}^{4}}+\EE{A_{k\Delta}^{4}+B_{k\Delta}^{4}+E_{k\Delta}^{4}}\right)\mathbb{P}\left(\Omega_{n}^{c}\right)\right]^{1/2}.\end{align*}
By Lemmas \ref{lem_majoration_AB} and \ref{majoration_omegaC}, we
obtain: \[
\E{\norm{\hat{g}_{m}-g_{A}}_{n}^{2}\units{\Omega_{n}^{c}}}\lesssim\frac{1}{\Delta^{3/2}n^{4}}\leq\frac{1}{n}.\]

\subsection{Proof of Theorem \ref{thme_est_sigma_xi_adaptatif}}

First, we apply the Berbee's coupling lemma to the random vectors
$(B_{k\Delta}+E_{k\Delta},X_{k\Delta})$ which are exponentially $\beta$-mixing.
According to Berbee's coupling lemma, we can construct independent
variables \[
U_{k,a}^{*}=\frac{1}{q_{n}}\sum_{l=1}^{q_{n}}\left(B+E\right)_{(2(k-1)+a)q_{n}+l)\Delta}^{*}t(X_{(2(k-1+a)q_{n}+l)\Delta}^{*})\]
 such that for $a\in\left\{ 0,1\right\} $, the random variables $(U_{k,a}^{*})_{0\leq k\leq p_{n}}$
are independent and have same law as \[
U_{k,a}=\frac{1}{q_{n}}\sum_{l=1}^{q_{n}}\left(B+E\right)_{(2(k-1)+a)q_{n}+l)\Delta}t(X_{(2(k-1+a)q_{n}+l)\Delta}).\]
 Let us set \[
\Omega^{*}=\left\{ \omega,\:\forall a,\:\forall k,\: U_{k,a}=U_{k,a}^{*}\right\} ,\]
 \[
\Omega_{B,\alpha}=\left\{ \omega,\:\forall(a,k),\:\left|U_{k,a}^{*}\right|\leq c\ln^{2}(n)D^{1/2}\Delta^{-\alpha}\right\} \]
 and \[
\rond{O}=\Omega_{n}\cap\Omega_{B,\alpha}\cap\Omega^{*}\]
 with $D=D_{m}+D_{m'}$. By Berbee's coupling lemma, \[
\mathbb{P}\left(\Omega^{*c}\right)\lesssim n\Delta/n^{8}.\]
 The following lemma is proved later.

\begin{lem}\label{lem_majoration_omega_B}

For any $\alpha>0$, there exists a constant $c$ such that  \[
\mathbb{P}\left(\Omega_{B,\alpha}^{c}\right)\lesssim\frac{1}{n^{5}}.\]

\end{lem}

Then \[
\mathbb{P}\left(\rond{O}^{c}\right)\leq\mathbb{P}(\Omega_{B,\alpha}^{c})+\mathbb{P}(\Omega^{*c})+\mathbb{P}(\Omega_{n}^{c})\lesssim\frac{1}{n^{5}}+\frac{n\Delta}{n^{8}}+\frac{1}{n^{8}}\lesssim\frac{1}{n^{5}}.\]
 We can bound $\E{\norm{\hat{g}_{\hat{m}}-g}_{n}^{2}\units{\rond{O}^{c}}}$
in the same way as we bound the risk of the non-adaptive estimator
on $\Omega_{n}^{c}$: \[
\E{\norm{\hat{g}_{\hat{m}}-g_{A}}_{n}^{2}\units{\rond{O}^{c}}}\lesssim\frac{1}{\Delta^{3/2}n^{5/2}}\lesssim\frac{1}{n}.\]
 It remains to bound the risk on $\rond{O}$. Let us set, for $a\in\left\{ 0,1\right\} $,
\[
\nu_{n,a}^{*}(t)=\frac{1}{p_{n}}\sum_{k=1}^{p_{n}}U_{k,a}^{*}\units{\rond{O}}-\E{U_{k,a}^{*}\units{\rond{O}}}\]
 and $\nu_{n}^{*}(t)=\nu_{n,0}^{*}(t)+\nu_{n,1}^{*}(t)$. We have:
\[
\norm{\hat{g}_{\hat{m}}-g}_{n}^{2}\units{\rond{O}}\lesssim\frac{1}{n}\sum_{k=1}^{n}A_{k\Delta}^{2}+\left(\EE{\left(B_{k\Delta}^{*}+E_{k\Delta}^{*}\right)\units{\rond{O}}}\right)^{2}+\sup_{t\in\rond{B}_{m,\hat{m}}}\left(\nu_{n}^{*}(t)\right)^{2}+2pen(m).\]
 As the random variables $B_{k\Delta}^{*}$ and $E_{k\Delta}^{*}$
are centred, \[
R_{1}:=\EE{\left(B_{k\Delta}^{*}+E_{k\Delta}^{*}\right)\units{\rond{O}}}=-\EE{\left(B_{k\Delta}^{*}+E_{k\Delta}^{*}\right)\units{\rond{O}^{c}}}\]
 then by Lemma \ref{lem_majoration_AB}, \[
\left|R_{1}\right|\lesssim\left(\EE{\left(B_{k\Delta}^{*}+E_{k\Delta}^{*}\right)^{2}}\mathbb{P}\left(\rond{O}^{c}\right)\right)^{1/2}\leq n^{-5/2}\Delta^{-1/2}.\]
 Then \begin{align*}
\E{\norm{\hat{g}_{\hat{m}}-g_{A}}_{n}^{2}\units{\rond{O}}} & \lesssim\norm{g_{m}-g}_{L^{2}}^{2}+\Delta+\frac{1}{n}+2pen(m)\\
 & +12\E{\sum_{m'\in\rond{M}_{n}}\left[\left(\sup_{t\in\rond{B}_{m,m'}}\left(\nu_{n}^{*}(t)\right)^{2}-p(m,m')\right)\units{\rond{O}}\right]_{+}}\end{align*}
 The functions $\nu_{n,a}^{*}(t)$ satisfy the assumptions of Talagrand's
inequality with $M=c\ln^{2}(n)D^{1/2}\Delta^{-\alpha}$, $V=\frac{1}{q_{n}\Delta}$,
and $H^{2}=c'D/(n\Delta)$. Then \begin{align*}
R_{2} & :=\E{\left[\sup_{t\in\rond{B}_{m,m'}}\left(\nu_{n}^{*}(t)\right)^{2}-12p(m,m')\right]_{+}}\\
 & \lesssim\frac{1}{n\Delta}\exp\left(-c'\frac{p_{n}q_{n}\Delta}{n\Delta}D\right)+\frac{c^{2}\ln^{4}(n)D\Delta^{-2\alpha}}{p_{n}^{2}}\exp\left(-\frac{c'}{c}\frac{p_{n}D^{1/2}}{\sqrt{n\Delta}\Delta^{-\alpha}D^{1/2}\ln^{2}(n)}\right)\\
 & \lesssim\frac{1}{n\Delta}\exp\left(-cD\right)+\frac{\ln^{6}(n)}{n^{2}\Delta^{2+2\alpha}}D\exp\left(-c'\frac{\sqrt{n\Delta^{1+2\alpha}}}{\ln^{3}(n)}\right).\end{align*}
Consequently, as $\alpha$ is as small as we want: \[
\E{\left[\sup_{t\in\rond{B}_{m,\hat{m}}}\nu_{n}^{2}(t)-pen(m)\right]_{+}}\lesssim\frac{1}{n\Delta}\sum_{m'}e^{-cD_{m,m'}}\lesssim\frac{1}{n\Delta}.\]

\subsection{Proof of Lemma \ref{lem_majoration_omega_B}}

We have that \[
U_{1,0}^{*}=\frac{1}{q_{n}}\sum_{j=1}^{q_{n}}(B_{j\Delta}^{*}+E_{j\Delta}^{*})t(X_{j\Delta}^{*})\lesssim\frac{1}{q_{n}\Delta}\sum_{j=1}^{q_{n}}\left(J_{j\Delta}^{2}+Z_{j\Delta}^{2}\right)t(X_{j\Delta}).\]
 We know that $\left|t(X_{j\Delta})\right|\leq\norm{t}_{\infty}\lesssim D^{1/2}$.
Moreover, \begin{align}
\mathbb{P}\left(\left|Z_{k\Delta}\right|\geq k\sigma_{0}\Delta^{1/2}\ln(n)\right) & \leq n^{-k}\EE{\exp\left(\frac{1}{\sigma_{0}\Delta^{1/2}}Z_{k\Delta}\right)}\nonumber \\
 & \leq n^{-k}\EE{\exp\left(\frac{1}{\sigma_{0}^{2}\Delta}\int_{k\Delta}^{(k+1)\Delta}\sigma^{2}(X_{s})ds\right)}\nonumber \\
 & \leq n^{-k}.\label{eq:majoration_exp_Zk}\end{align}
Then \begin{equation}
\mathbb{P}\left(\left|Z_{k\Delta}\right|\geq6\sigma_{0}^{2}\Delta\ln(n)\right)\lesssim n^{-6}.\label{eq:majoration_Z}\end{equation}
 and then \[
\sum_{k=1}^{p_{n}}\mathbb{P}\left(\frac{1}{q_{n}}\sum_{j=1}^{q_{n}}Z_{j\Delta}^{2}\geq36\sigma_{0}^{4}\Delta\ln^{2}(n)\right)\lesssim n^{-5}.\]

\subsubsection*{Bound of $\mathbb{P}\left(\left|J_{k\Delta}^{(1)}\right|\geq12\xi_{0}\Delta^{1/2}\ln(n)\right)$.}

The terms $J_{k\Delta}^{(1)}$ are small and can be bounded in the
same way as the Brownian terms $Z_{k\Delta}$. As $\nu$ is symmetric:
\begin{align*}
\mathbb{P}\left(\left|J_{k\Delta}^{(1)}\right|\geq12\xi_{0}\Delta^{1/2}\ln(n)\right) & \leq2\mathbb{P}\left(\exp\left(aJ_{k\Delta}^{(1)}\right)\geq\exp\left(12a\xi_{0}\Delta^{1/2}\ln(n)\right)\right)\\
 & \leq2\exp\left(-12a\xi_{0}\Delta^{1/2}\ln(n)\right)\E{\exp\left(aJ_{k\Delta}^{(1)}\right)}.\end{align*}
 According to Corollary 5.2.2 of \citet{applebaum}, \[
\E{\exp\left(aJ_{k\Delta}^{(1)}\right)}=\E{\exp\left(\int_{k\Delta}^{(k+1)\Delta}\int_{-\Delta^{1/2}}^{\Delta^{1/2}}\left(e^{a\xi(X_{s-})z}-1-a\xi(X_{s^{-}})z\right)\nu(dz)ds\right)}.\]
 Then for any $a\leq1/(2\xi_{0}\Delta^{1/2})$,

\begin{align*}
\E{\exp\left(aJ_{k\Delta}^{(1)}\right)} & \leq\E{\exp\left(\int_{k\Delta}^{(k+1)\Delta}\int_{-\Delta^{1/2}}^{\Delta^{1/2}}a^{2}z^{2}\xi^{2}(X_{s^{-}})\nu(dz)\right)}\\
 & \leq\E{\exp\left(\xi_{0}^{2}a^{2}\Delta^{\text{2-\ensuremath{\beta}/2}}\right)}.\end{align*}
 Let us then set $a=1/(2\xi_{0}\Delta^{1/2})$, we obtain: \begin{equation}
\mathbb{P}\left(\left|J_{k\Delta}^{(1)}\right|\geq12\xi_{0}\Delta^{1/2}\ln(n)\right)\lesssim\exp\left(-6\ln(n)\right)\leq n^{-6}.\label{eq:majoration_j1}\end{equation}

\subsubsection*{Bound for the jumps greater than $\Delta^{1/2}$.}

The probability that \[
J_{k\Delta}^{(2)}+J_{k\Delta}^{(3)}\geq\Delta^{1/2}\ln(n)\]
 is not small enough. We have to bound both the number of jumps of
the time interval $[k\Delta,(k+1)\Delta[$ and the size of the jumps.
Let us first consider the jumps greater than 1: \[
J_{k\Delta}^{(0)}=\int_{k\Delta}^{(k+1)\Delta}\xi(X_{s^{-}})\int_{\left|z\right|\geq1}\mu(dz,ds).\]
 The probability of having a very high jump is quite small:~by Assumption
A\ref{hypo_nu_sous_exp}, \begin{equation}
\nu\left(\left[-\frac{8\ln(n)}{\lambda},\frac{8\ln(n)}{\lambda}\right]^{c}\right)\lesssim n^{-8}.\label{eq:majoration_taille_sauts}\end{equation}
 The probability of having more than $C=8\eta/(1-\beta/2)$ (see Assumption
A\ref{hypo_nu_sous_exp}) jumps greater than $1$ on a time interval
$\Delta$ is very low: \begin{align}
Q_{1} & :=\mathbb{P}\left(\mu\left(\left[k\Delta,(k+1)\Delta\right[,\left[-1,1\right]^{c}\right)\geq C\right)\nonumber \\
 & \leq\mathbb{P}\left(\mu\left(\left[k\Delta,(k+1)\Delta\right[,\left[-\Delta^{1/2},\Delta^{1/2}\right]^{c}\right)\geq C\right)\nonumber \\
 & \leq\left(\Delta\int_{\left|z\right|>\Delta^{1/2}}\nu(dz)\right)^{C}\lesssim\left(\Delta^{1-\beta/2}\right)^{C}=\Delta^{8\eta}\lesssim n^{-8}.\label{eq:majoration_nb_sauts_Delta}\end{align}
 By \eqref{eq:majoration_taille_sauts} and \eqref{eq:majoration_nb_sauts_Delta},
\[
\mathbb{P}\left(\vert J_{k\Delta}^{(0)}\vert\geq\frac{8C\ln(n)}{\lambda}\right)\lesssim n^{-8}.\]
 Let us set $v_{1}=\nu(]-1,1[^{c})\vee1$. We have that \begin{align}
Q_{2} & :=\mathbb{P}\left(\frac{1}{q_{n}}\sum_{k=1}^{q_{n}}\left|J_{k\Delta}^{(0)}\right|^{2}\geq\frac{8^{2}C^{2}}{\lambda^{2}}v_{1}\Delta\ln^{2}(n)\right)\nonumber \\
 & \lesssim q_{n}\mathbb{P}\left(\vert J_{k\Delta}^{(0)}\vert\geq\frac{8C\ln(n)}{\lambda}\right)+\mathbb{P}\left[\mu\left([0,q_{n}\Delta[,[-1,1]^{c}\right)\geq v_{1}\Delta q_{n}\right]\nonumber \\
 & \lesssim\frac{q_{n}}{n^{8}}+\sum_{j\geq v_{1}\Delta q_{n}}^{+\infty}\frac{\left(q_{n}\Delta v_{1}\right)^{j}}{j!}e^{-q_{n}\Delta v_{1}}\nonumber \\
 & \lesssim\frac{q_{n}}{n^{8}}+\sum_{j\geq8v_{1}\ln(n)}^{+\infty}\left(\frac{8\ln(n)v_{1}e}{j}\right)^{j}\sqrt{j}e^{-8\ln(n)v_{1}}\nonumber \\
 & \lesssim\frac{q_{n}}{n^{8}}+\frac{1}{n^{8}}.\label{eq:majoration_J0}\end{align}
 Let us now set $\alpha_{0}=0$, $\alpha_{j}=\frac{2\alpha_{j-1}+\alpha}{\beta}\wedge\frac{1}{2}$
and \[
J_{k\Delta}^{(\alpha_{j})}=\int_{k\Delta}^{(k+1)\Delta}\int_{[-\Delta^{\alpha_{j-1}},-\Delta^{\alpha_{j}}]\cup[\Delta^{\alpha_{j}},\Delta^{\alpha_{j-1}}]}\xi(X_{s^{-}})dL_{s}.\]
 By \eqref{eq:majoration_nb_sauts_Delta}, \[
\mathbb{P}\left(\left|J_{k\Delta}^{(\alpha_{j})}\right|\geq C\Delta^{\alpha_{j-1}}\right)\lesssim\frac{1}{n^{8}}.\]
 We have that $\nu\left(\left[-\Delta^{-\alpha_{j-1}},-\Delta^{\alpha_{j}}\right]\cup\left[\Delta^{\alpha_{j}},\Delta^{\alpha_{j-1}}\right]\right)\lesssim\Delta^{-\alpha_{j}\beta}$.
Let us set $v_{2}=\Delta^{\beta\alpha_{j}}\left(\nu(\left[-\Delta^{-\alpha_{j-1}},-\Delta^{\alpha_{j}}\right]\cup\left[\Delta^{\alpha_{j}},\Delta^{\alpha_{j-1}}\right])\vee1\right).$
Then \begin{align}
Q_{3} & :=\mathbb{P}\left(\frac{1}{q_{n}}\sum_{k=1}^{q_{n}}\left|J_{k\Delta}^{\left(\alpha_{j}\right)}\right|^{2}\geq C^{2}\Delta^{1-\alpha}\right)\nonumber \\
 & \leq\mathbb{P}\left[\mu\left(\left[0,q_{n}\Delta\right[,\left[-\Delta^{-\alpha_{j-1}},-\Delta^{\alpha_{j}}\right]\cup\left[\Delta^{\alpha_{j}},\Delta^{\alpha_{j-1}}\right]\right)\geq q_{n}\Delta^{1-2\alpha_{j-1}-\alpha}\right]\nonumber \\
 & +q_{n}\mathbb{P}\left(\left|J_{k\Delta}^{(\alpha_{j})}\right|\geq C\Delta^{\alpha_{j-1}}\right)\nonumber \\
 & \lesssim\frac{q_{n}}{n^{8}}+\sum_{i=q_{n}\Delta^{1-2\alpha_{j-1}-\alpha}}^{\infty}\frac{\left(q_{n}\Delta v_{2}\Delta^{-\beta\alpha_{j}}\right)^{i}}{i!}e^{-q_{n}\Delta v_{2}\Delta^{-\beta\alpha_{j}}}\nonumber \\
 & \lesssim n^{-8}.\label{eq:majoration_Jj}\end{align}
 Then, by \eqref{eq:majoration_Z}, \eqref{eq:majoration_j1}, \eqref{eq:majoration_J0}
and \eqref{eq:majoration_Jj}, we obtain: \[
\mathbb{P}\left(\Omega_{B,\alpha}^{c}\right)\lesssim\frac{p_{n}q_{n}}{n^{6}}+\frac{p_{n}q_{n}}{n^{8}}\lesssim n^{-5}.\]

\subsection{Proof of Lemma \ref{lem_majoration_omega_X}}

\paragraph{Bound of $\mathbb{P}\left(\Omega_{X,k}^{c}\right)$. }

We have that $X_{(k+1)\Delta}=X_{k\Delta}+\int_{k\Delta}^{(k+1)\Delta}b(X_{s})ds+Z_{k\Delta}+J_{k\Delta}$.
Then

\begin{align*}
\mathbb{P}\left(\Omega_{X,k}^{c}\right) & \leq\mathbb{P}\left(\left|\int_{k\Delta}^{(k+1)\Delta}b(X_{s})ds\right|\geq\Delta^{1/2}\right)+\mathbb{P}\left(\left|Z_{k\Delta}\right|\geq\sigma_{0}\Delta^{1/2}\ln(n)\right)\\
 & +\mathbb{P}\left(\left|J_{k\Delta}\right|\geq\xi_{0}\Delta^{1/2}\ln(n)\right).\end{align*}
 By Markov's inequalities, for any $k\leq4$: \begin{equation}
\mathbb{P}\left(\left|\int_{k\Delta}^{(k+1)\Delta}b(X_{s})ds\right|\geq\Delta^{1/2}\right)\lesssim\Delta^{-k}\EE{\left(\int_{k\Delta}^{(k+1)\Delta}b(X_{s})ds\right)^{2k}}\lesssim\Delta^{k}\label{eq:majoration_b}\end{equation}
 and by \eqref{eq:majoration_exp_Zk}, $\mathbb{P}\left(\left|Z_{k\Delta}\right|\geq k\sigma_{0\mbox{}}\Delta^{1/2}\ln(n)\right)\lesssim n^{-k}.$
 Moreover, \begin{align}
\mathbb{P}\left(\left|J_{k\Delta}^{(2)}+J_{k\Delta}^{(3)}\right|>0\right) & \leq\Delta\int_{[-\Delta^{1/2},\Delta^{1/2}]^{c}}\nu(dz)\leq\Delta\Delta^{-\beta/2}\int_{[-\Delta^{1/2},\Delta^{1/2}]^{c}}z^{\beta}\nu(dz)\nonumber \\
 & \lesssim\Delta^{1-\beta/2}\label{eq:majoration_Jk2_Jk3}\end{align}
 and by Markov's inequality: \begin{align}
\mathbb{P}\left(\left|J_{k\Delta}^{(1)}\right|>\xi_{0}\Delta^{1/2}\right) & \leq\frac{1}{\xi_{0}^{2}\Delta}\EE{\left(J_{k\Delta}^{(1)}\right)^{2}}\nonumber \\
 & \leq\frac{1}{\xi_{0}^{2}\Delta}\Delta\xi_{0}^{2}\int_{-\Delta^{1/2}}^{\Delta^{1/2}}z^{2}\nu(dz)\lesssim\Delta^{1-\beta/2}\int_{-\Delta^{1/2}}^{\Delta^{1/2}}z^{\beta}\nu(dz)\nonumber \\
 & \lesssim\Delta^{1-\beta/2}.\label{eq:majoration_Jk1}\end{align}

\subsubsection*{Bound of $\mathbb{P}\left(\Omega_{N,k}^{c}\right)$. }

We have that \[
\mathbb{P}\left(N_{k}\geq1\right)=\int_{k\Delta}^{(k+1)\Delta}\int_{\vert z\vert\geq\Delta^{1/4}}\nu(dz)\leq\Delta^{1-\beta/4}\int_{\vert z\vert\geq\Delta^{1/4}}z^{\beta}\nu(dz)\lesssim\Delta^{1-\beta/4}.\]
 Then by \eqref{eq:majoration_Jk2_Jk3} and \eqref{eq:majoration_Jk1},
we obtain: \begin{align*}
\mathbb{P}\left(\Omega_{N,k}^{c}\right) & \leq\mathbb{P}\left(N_{k}\geq1\right)+\mathbb{P}\left(\vert J_{k\Delta}^{(2)}+J_{k\Delta}^{(3)}\vert>0\right)+\mathbb{P}\left(\left|\int_{k\Delta}^{(k+1)\Delta}dL_{s}^{(1)}\right|\geq\ln(n)\Delta^{1/2}\right)\\
 & \lesssim\Delta^{1-\beta/2}.\end{align*}

\subsubsection*{Bound of $\mathbb{P}\left(\Omega_{X,k}\cap\Omega_{N,k}^{c}\right)$. }

We have that \[
\mathbb{P}\left(\Omega_{X,k}\cap\left\{ N_{k}\geq1\right\} \right)\leq\mathbb{P}\left(N_{k}\geq2\right)+\mathbb{P}\left(\Omega_{X,k}\cap\left\{ N_{k}=1\right\} \right).\]
 Now $\mathbb{P}\left(N_{k}\geq2\right)\leq\left(\Delta^{1-\beta/4}\int_{\vert z\vert\geq\Delta^{1/4}}z^{\beta}\nu(dz)\right)^{2}\lesssim\Delta^{2-\beta/2}$.
Moreover, if $N_{k}=1$, then $\left|J_{k\Delta}^{(3)}\right|\geq\xi_{1}\Delta^{1/4}$
and by conditional independence, we get: \begin{align*}
S_{1} & :=\mathbb{P}\left(\Omega_{X,k}\cap\left\{ N_{k}=1\right\} \right)\\
 & \leq\mathbb{P}\left(N_{k}=1\right)\times\mathbb{P}\left(\left|\int_{k\Delta}^{(k+1)\Delta}b(X_{s})ds+Z_{k\Delta}+J_{k\Delta}^{(1)}+J_{k\Delta}^{(2)}\right|>\xi_{1}\Delta^{1/4}\right)\\
 & \leq\mathbb{P}\left(N_{k}=1\right)\left[\mathbb{P}\left(\left|\int_{k\Delta}^{(k+1)\Delta}b(X_{s})ds\right|\geq\frac{\xi_{1}\Delta^{1/4}}{3}\right)+\mathbb{P}\left(\left|Z_{k\Delta}\right|\geq\frac{\xi_{1}\Delta^{1/4}}{3}\right)\right]\\
 & +\mathbb{P}\left(N_{k}=1\right)\mathbb{P}\left(\left|J_{k\Delta}^{(1)}+J_{k\Delta}^{(2)}\right|\geq\frac{\xi_{1}\Delta^{1/4}}{3}\right).\end{align*}
 By \eqref{eq:majoration_b} and \eqref{eq:majoration_exp_Zk}, $\mathbb{P}\left(\left|\int_{k\Delta}^{(k+1)\Delta}b(X_{s})ds\right|\geq c\Delta^{1/4}\right)\leq\Delta^{4}$
and, as $\ln(n)\ll\Delta^{-1/4}$, $\mathbb{P}\left(\left|Z_{k\Delta}\right|\geq c\Delta^{1/4}\right)\leq\mathbb{P}\left(\left|Z_{k\Delta}\right|\geq\ln(n)\Delta^{1/2}\right)\lesssim n^{-1}$.
Moreover, by a Markov inequality, we obtain: \begin{align*}
\mathbb{P}\left(\left|J_{k\Delta}^{(1)}+J_{k\Delta}^{(2)}\right|>c\Delta^{1/4}\right) & \leq c^{-2}\Delta^{-1/2}\EE{\left(J_{k\Delta}^{(1)}+J_{k\Delta}^{(2)}\right)^{2}}\\
 & \leq c^{-2}\Delta^{-1/2}\xi_{0}^{2}\Delta\int_{-\Delta^{1/4}}^{\Delta^{1/4}}z^{2}\nu(dz)\\
 & \lesssim\Delta^{1/2}\Delta^{1/2-\beta/4}.\end{align*}
 As $\mathbb{P}\left(N_{k}=1\right)\lesssim\Delta^{1-\beta/4}$, we
obtain:~ \[
\mathbb{P}\left(\Omega_{X,k}\cap\left\{ N_{k}\geq1\right\} \right)\lesssim\Delta^{2-\beta/2}.\]
 Let us set $L_{s}^{(1)+(2)}=L_{s}^{(1)}+L_{s}^{(2)}$ and $J_{k\Delta}^{(1)+(2)}=J_{k\Delta}^{(1)}+J_{k\Delta}^{(2)}$.
We consider \[
\rond{E}_{k}=\left\{ \left|\int_{k\Delta}^{(k+1)\Delta}dL_{s}^{(1)+(2)}\right|\leq4\frac{\xi_{0}+\sigma_{0}}{\xi_{1}}\Delta^{1/2}\ln(n)\right\} .\]
 We have that \begin{align*}
\rond{E}_{k}^{c} & \subseteq\left\{ \left|\xi(X_{k\Delta^{-}})\int_{k\Delta}^{(k+1)\Delta}dL_{s}^{(1)+(2)}\right|\geq4(\xi_{0}+\sigma_{0})\Delta^{1/2}\ln(n)\right\} \\
 & \subseteq\left\{ \left|J_{k\Delta}^{(1)+(2)}\right|\geq2\left(\xi_{0}+\sigma_{0}\right)\Delta^{1/2}\ln(n)\right\} \\
 & \cup\left\{ \left|\int_{k\Delta}^{(k+1)\Delta}\left(\xi(X_{s^{-}})-\xi(X_{k\Delta^{-}})\right)dL_{s}^{(1)+(2)}\right|\geq2\left(\xi_{0}+\sigma_{0}\right)\Delta^{1/2}\ln(n)\right\} .\end{align*}
 By \eqref{eq:majoration_b} and \eqref{eq:majoration_exp_Zk}, \begin{align*}
S_{2} & :=\mathbb{P}\left(\Omega_{X,k}\cap\left\{ \left|J_{k\Delta}^{(1)+(2)}\right|\geq2\left(\sigma_{0}+\xi_{0}\right)\Delta^{1/2}\ln(n)\right\} \cap N_{k}=0\right)\\
 & \lesssim\mathbb{P}\left(\left|\int_{k\Delta}^{(k+1)\Delta}b(X_{s})ds+Z_{k\Delta}\right|\geq\left(\sigma_{0}+\xi_{0}\right)\Delta^{1/2}\ln(n)\right)\\
 & \lesssim\Delta^{4}+n^{-1}.\end{align*}
 By the \burk, we obtain that \[
\E{\sup_{s\leq(k+1)\Delta}\left(X_{s}-X_{k\Delta}\right)^{4}\units{\Omega_{N,k}}}\lesssim\Delta^{2-\beta/4}.\]
 Moreover, \begin{align}
S_{3} & :=\EE{\left(\int_{k\Delta}^{(k+1)\Delta}\left(\xi(X_{s^{-}})-\xi(X_{k\Delta^{-}})\right)dL_{s}^{(1)+(2)}\right)^{4}}\nonumber \\
 & \lesssim\left(\int_{k\Delta}^{(k+1)\Delta}\Delta\int_{-\Delta^{1/4}}^{\Delta^{1/4}}z^{2}\nu(dz)\right)^{2}+\int_{k\Delta}^{(k+1)\Delta}\Delta^{2-\beta/4}\int_{-\Delta^{1/4}}^{\Delta^{1/4}}z^{4}\nu(dz)\nonumber \\
 & \lesssim\Delta^{5-\beta/2}+\Delta^{4-\beta/2}\lesssim\Delta^{4-\beta/2}.\label{eq:majoration_xi_4}\end{align}
 Then by a Markov's inequality, \[
\mathbb{P}\left(\left|\int_{k\Delta}^{(k+1)\Delta}\left(\xi(X_{s^{-}})-\xi(X_{k\Delta^{-}})\right)dL_{s}^{(1)+(2)}\right|\geq\Delta^{1/2}\ln(n)\right)\lesssim\Delta^{2-\beta/2}\]
 which ends the proof.

\subsection{Proof of Lemma \ref{lem_majoration_AB-1}}

From the \burk and Proposition \ref{prop_Xk}, we derive easily the
bounds for $\tilde{A}_{k\Delta}$ and $\tilde{B}_{k\Delta}$. It remains
to bound $\E{\tilde{E}_{k\Delta}^{2}\left|\rond{F}_{k\Delta}\right.}$
and $\E{\tilde{E}_{k\Delta}^{4}\left|\rond{F}_{k\Delta}\right.}$.
We first bound $\E{\left(J_{k\Delta}^{(1)+(2)}\right)^{4}}$. We have
that \[
J_{k\Delta}^{(1)+(2)}=\int_{k\Delta}^{(k+1)\Delta}(\xi(X_{s^{-}})-\xi(X_{k\Delta^{-}}))dL_{s}^{(1)+(2)}+\xi(X_{k\Delta^{-}})\int_{k\Delta}^{(k+1)\Delta}dL_{s}^{(1)+(2)}.\]
 By \eqref{eq:majoration_xi_4}, $\EE{\left(\int_{k\Delta}^{(k+1)\Delta}(\xi(X_{s^{-}})-\xi(X_{k\Delta^{-}}))dL_{s}^{(1)+(2)}\right)^{4}}\lesssim\Delta^{4-\beta/2}$.
It remains to bound $\EE{\left(\int_{k\Delta}^{(k+1)\Delta}dL_{s}^{(1)+(2)}\units{\rond{E}_{k}}\right)^{4}}.$
This is nearly Proposition 4.5 of \citet{mai2012}. Let us introduce
a nonnegative function $f$ $\rond{C}^{\infty}$ such that \[
\begin{cases}
f(x)=x^{4} & \textrm{ if }\vert x\vert\leq1\\
f(x)=0 & \textrm{ if }\vert x\vert\geq2.\end{cases}\]
 Let us set $f^{a}(x)=a^{4}f(x/a)$. By stationarity, we have \[
\EE{\left(\int_{k\Delta}^{(k+1)\Delta}dL_{s}^{(1)+(2)}\units{\Omega_{N,k}}\right)^{4}}=\EE{\left(L_{\Delta}^{(1)+(2)}\right)^{4}\units{\Omega_{N,k}}}\leq\E{f^{\Delta^{1/2}\ln(n)}(L_{\Delta}^{(1)+(2)})}.\]
The following result is needed. 

\begin{res}{[}Fourier transform{]}\label{res_fourier}

We denote by $\mathcal{F}h$ the Fourier transform of a function $h\in L^{1}(\mathbb{R})$:
\[
\mathcal{F}h(x)=\int_{\mathbb{R}}f(u)e^{-ixu}du.\]
 The Schwarz space is defined as \[
S\left(\mathbb{R}\right)=\left\{ h\in\rond{C}^{\infty},\;\forall p,q\in\mathbb{N},\;\exists C_{p,q},\:\forall x\in\mathbb{R},\;\vert x^{p}h^{(q)}(x)\vert\leq C_{pq}\right\} .\]

Then we have the following properties:~ 
\begin{enumerate}
\item For any $h_{1},h_{2}\in L^{2}(\mathbb{R})$, $(a_{1},a_{2})\in\mathbb{R}^{2}$,
$\mathcal{F}(a_{1}h_{1}+a_{2}h_{2})=a_{1}\mathcal{F}h_{1}+a_{2}\mathcal{F}h_{2}$.\label{enu:fourier_linearite} 
\item For any $h\in L^{2}(\mathbb{R})$, $\mathcal{F}h\in L^{2}(\mathbb{R})$
and $\forall x\in\mathbb{R},\; h(x)=\frac{1}{2\pi}\int_{\mathbb{R}}e^{itx}\mathcal{F}h(t)dt$.
\label{enu:fourier_inverse} 
\item For any $h\in L^{2}(\mathbb{R})$, $\mathcal{F}h(./a)(x)=\vert a\vert\mathcal{F}h(ax)$.
\label{enu:fourier_dilatation} 
\item For any functions $h_{1},h_{2}\in L^{2}(\mathbb{R})$, the Parseval's
formula holds:\label{enu:fourier_parseval} \[
\int_{\mathbb{R}}h_{1}(x)\overline{h_{2}(x)}dx=\frac{1}{2\pi}\int_{\mathbb{R}}\mathcal{F}h_{1}(u)\overline{\mathcal{F}h_{2}}(u)du.\]
 As $\mathcal{F}\delta_{y}(x)=e^{-ixy}$, \[
h(0)=\int_{\mathbb{R}}h(y)\delta_{0}(y)dy=\frac{1}{2\pi}\int_{\mathbb{R}}\mathcal{F}h(u)du\]

\item For any $h$ in $S(\mathbb{R})$, $\mathcal{F}h\in S(\mathbb{R})$
and \label{enu:fourier_derivee} \begin{gather*}
\mathcal{F}(h^{(q)})(x)=(ix)^{q}\mathcal{F}h(x).\end{gather*}

\end{enumerate}
\end{res}

By Result \ref{res_fourier}.\ref{enu:fourier_parseval}, we have
that \begin{align*}
\E{f^{a}(L_{t}^{(1)+(2)})} & =\int_{\mathbb{R}}f^{a}(x)P_{L_{\Delta}^{(1)+(2)}}(dx)\\
 & =\frac{1}{2\pi}\int_{\mathbb{R}}\mathcal{F}f^{a}(u)\bar{\phi}_{\Delta}(u)du\end{align*}
 where $\phi_{\Delta}$ is the characteristic function of the Lévy
process $L_{\Delta}^{(1)+(2)}$: \[
\phi_{\Delta}(u)=\exp\left(\Delta\int_{-\Delta^{1/4}}^{\Delta^{1/4}}(e^{iux}-1-iux)\nu(dx)\right).\]
 By a Taylor development in 0, we obtain that \[
\phi_{\Delta}(u)=1+\psi_{\Delta}(u)+R(\Delta,u)\]
 with $\psi_{\Delta}(u)=\Delta\int_{-\Delta^{1/4}}^{\Delta^{1/4}}(e^{iux}-1-iux)\nu(dx)$.
Then \begin{alignat*}{1}
\E{f^{a}(L_{t}^{(1)+(2)})} & =\frac{1}{2\pi}\int_{\mathbb{R}}\mathcal{F}f^{a}(u)du+\frac{1}{2\pi}\int_{\mathbb{R}}\mathcal{F}f^{a}(u)\overline{\psi_{\Delta}(u)}du\\
 & +\frac{1}{2\pi}\int_{\mathbb{R}}\mathcal{F}f^{a}(u)\overline{R(\Delta,u)}du.\end{alignat*}
 By Result \ref{res_fourier}.\ref{enu:fourier_parseval}, $\int_{\mathbb{R}}\mathcal{F}f^{a}(u)du=2\pi f^{a}(0)=0$
and consequently, \begin{align*}
\int_{\mathbb{R}}\mathcal{F}f^{a}(u)\overline{\psi_{\Delta}(u)}du & =\int_{\mathbb{R}}\mathcal{F}f^{a}(u)\Delta\int_{-\Delta^{1/4}}^{\Delta^{1/4}}(e^{-iux}-1+iux)du\nu(dx)\\
 & =\Delta\int_{-\Delta^{1/4}}^{\Delta^{1/4}}(2\pi)f^{a}(x)\nu(dx)-f^{a}(0)+\int_{\mathbb{R}}\mathcal{F}f^{a}(u)iuxdu\nu(dx).\end{align*}
 By Result \ref{res_fourier}.\ref{enu:fourier_derivee}, as $f^{a}\in S(\mathbb{R})$,
$\int_{\mathbb{R}}\mathcal{F}f^{a}(u)iudu=\int_{\mathbb{R}}\mathcal{F}((f^{a})')(u)du=(f^{a})'(0)=0$.
Then \begin{align*}
\int_{\mathbb{R}}\mathcal{F}f^{\Delta^{1/2}\ln(n)}(u)\overline{\psi_{\Delta}(u)}du & =2\pi\Delta\int_{-\Delta^{1/4}}^{\Delta^{1/4}}f^{\Delta^{1/2}\ln(n)}(x)\nu(dx)\\
 & \lesssim\Delta\int_{-2\Delta^{1/2}\ln(n)}^{2\Delta^{1/2}\ln(n)}x^{4}\nu(dx)\lesssim\Delta\Delta^{2-\beta/2}\ln(n)^{4-\beta}.\end{align*}
 It remains to bound $\E{\int_{\mathbb{R}}\mathcal{F}f^{a}(u)\overline{R(\Delta,u)}du}$.
We have that \[
\vert R(\Delta,u)\vert=\left|e^{\psi_{\Delta}(u)}-\psi_{\Delta}(u)-1\right|\leq\left|\psi_{\Delta}^{2}(u)\right|.\]
 According to \citet{kappus_these}, $\left|\psi_{\Delta}(u)\right|\lesssim C\Delta\vert u\vert^{\beta}$.
By Result \ref{res_fourier}.\ref{enu:fourier_dilatation}, $\mathcal{F}f^{a}(u)=a^{5}\mathcal{F}f(au)$
and therefore \begin{align*}
\left|\E{\int_{\mathbb{R}}\mathcal{F}f^{a}(u)\overline{R(\Delta,u)}du}\right| & \lesssim\Delta^{2}\int_{\mathbb{R}}\left|\mathcal{F}f^{a}(u)\right|\vert u\vert^{2\beta}du\\
 & \lesssim\Delta^{2}\int_{\mathbb{R}}a^{5}\left|\mathcal{F}f(au)\right|\vert u\vert^{2\beta}du.\end{align*}
 As $f^{a}\in S(\mathbb{R})$, $\mathcal{F}f^{a}\in S(\mathbb{R})$
and then for any $m>0$, $\exists C_{m}>0$, $\vert\mathcal{F}f(u)\vert\leq C_{m}\vert u\vert^{-m}$.
Then, for any $m\in\mathbb{N}$:~ \[
\E{\int_{\mathbb{R}}\mathcal{F}f^{\Delta^{1/2}\ln(n)}(u)\overline{R(\Delta,u)}du}\lesssim\Delta^{2}\int_{\mathbb{R}}a^{5-m}\vert u\vert^{2\beta-m}\wedge a^{5}\vert u\vert^{2\beta}du.\]
 We choose $m$ such that $2\beta+1<m\leq3+\beta$. As $\beta<2$,
$m$ always exists. Then $\int_{\mathbb{R}}\vert u\vert^{2\beta-m}\wedge\vert u\vert^{2\beta}<\infty$
and we get: \[
\E{\int_{\mathbb{R}}\mathcal{F}f^{\Delta^{1/2}\ln(n)}(u)\overline{R(\Delta,u)}du}\lesssim\Delta^{2}\left(\Delta^{1/2}\ln(n)\right)^{2-\beta}\lesssim\Delta^{3-\beta/2}\ln(n)^{2-\beta}\]
 Then we obtain \begin{align}
\EE{\left(J_{k\Delta}^{(1)+(2)}\right){}^{4}\units{\rond{E}_{k}}} & \lesssim\Delta^{3-\beta/2}\left(\ln(n)\right)^{4-\beta}.\label{eq:majoration_J_4}\end{align}

\subsubsection*{Bound of $\E{\tilde{E}_{k\Delta}^{2}\units{\Omega_{X,k}\cap\Omega_{N,k}}}$. }

On $\Omega_{N,k}$, \[
\tilde{E}_{k\Delta}=\left(2b(X_{k\Delta})J_{k\Delta}^{(1)+(2)}+\frac{J_{k\Delta}^{(1)+(2)}}{\Delta}Z_{k\Delta}+\frac{\left(J_{k\Delta}^{(1)+(2)}\right)^{2}}{\Delta}\right)\units{\rond{E}_{k}}.\]
 Then by \eqref{eq:majoration_J_4}, \begin{align*}
\E{\tilde{E}_{k\Delta}^{2}\units{\Omega_{N,k}\cap\Omega_{X,k}}} & \lesssim\frac{\E{(J_{k\Delta}^{(1)+(2)})^{4}\units{\rond{E}_{k}}}}{\Delta^{2}}+\E{Z_{k\Delta}^{4}}\\
 & \lesssim\Delta^{1-\beta/2}\ln(n)^{4-\beta}.\end{align*}

\subsubsection*{Bound of $\E{\tilde{E}_{k\Delta}^{2}\units{\Omega_{X,k}\cap\Omega_{N,k}^{c}}}$. }

We have that \begin{align*}
S_{5} & :=\E{\tilde{E}_{k\Delta}^{2}\units{\Omega_{X,k}\cap\Omega_{N,k}^{c}}}\\
 & \leq\E{\tilde{E}_{k\Delta}^{2}\units{\Omega_{X,k}\cap\Omega_{N,k}^{c}\cap\left\{ \left|J_{k\Delta}\right|+\left|Z_{k\Delta}\right|\geq9\sigma_{0}\ln(n)\Delta^{1/2}\right\} }}\\
 & +\E{\tilde{E}_{k\Delta}^{2}\units{\Omega_{X,k}\cap\Omega_{N,k}^{c}\cap\left\{ \left|J_{k\Delta}\right|+\left|Z_{k\Delta}\right|\leq9\sigma_{0}\ln(n)\Delta^{1/2}\right\} }}\\
 & \lesssim\left(\E{\tilde{E}_{k\Delta}^{4}}\mathbb{P}\left(\Omega_{X,k}\cap\left\{ \left|J_{k\Delta}\right|+\left|Z_{k\Delta}\right|\geq9\sigma_{0}\ln(n)\Delta^{1/2}\right\} \right)\right)^{1/2}\\
 & +\ln^{2}(n)\mathbb{P}\left(\Omega_{X,k}\cap\Omega_{N,k}^{c}\right).\end{align*}
 It remains to bound $\mathbb{P}\left(\Omega_{X,k}\cap\left\{ \left|J_{k\Delta}\right|+\left|Z_{k\Delta}\right|\geq9\sigma_{0}\ln(n)\Delta^{1/2}\right\} \right)$
. By inequality \eqref{eq:majoration_exp_Zk}, $\mathbb{P}\left(\left|Z_{k\Delta}\right|\geq4\sigma_{0}\ln(n)\Delta^{1/2}\right)\lesssim n^{-4}$
and \begin{align*}
S_{6} & :=\mathbb{P}\left(\left\{ \left|J_{k\Delta}\right|\geq5\sigma_{0}\Delta^{1/2}\ln(n)\right\} \cap\Omega_{X,k}\right)\\
 & \leq\mathbb{P}\left(\left|\int_{k\Delta}^{(k+1)\Delta}b(X_{s})ds+Z_{k\Delta}\right|\geq4\sigma_{0}\Delta^{1/2}\ln(n)\right)\\
 & \lesssim\Delta^{5}+n^{-4}.\end{align*}
 It follows that \[
\E{\tilde{E}_{k\Delta}^{2}\units{\Omega_{X,k}\cap\Omega_{N,k}^{c}}}\lesssim\left(\left(n^{-4}+\Delta^{5}\right)\E{\tilde{E}_{k\Delta}^{4}}\right)^{1/2}+\Delta^{2-\beta/2}\ln^{4}(n).\]
 As $\E{\tilde{E}_{k\Delta}^{4}}\lesssim1/\Delta^{3}$, we get that
$\E{\tilde{E}_{k\Delta}^{2}\units{\Omega_{X,k}\cap\Omega_{N,k}^{c}}}\lesssim\Delta^{1-\beta/2}\ln^{4}(n)$.

\subsection{Proof of Theorem \ref{thme_sigma_m_fixe}}

As before, we decompose the bound of the risk on $\Omega_{n}$ and
$\Omega_{n}^{c}$. We bound the risk on $\Omega_{n}^{c}$ in the same
way as in the proof of Theorem \ref{thme_risque_est_sigma_xi_m_fixe}.
On $\Omega_{n}$ , we obtain that:~ \[
\E{\norm{\hat{\sigma}_{m}^{2}-\sigma_{A}^{2}}_{n}^{2}\units{\Omega_{n}}}\leq3\norm{\sigma_{m}^{2}-\sigma_{A}^{2}}_{\pi}^{2}+12\E{F_{\Delta}^{2}}+12\E{\sup_{t\in\rond{B}_{m}}\tilde{\nu}_{n}^{2}(t)}.\]
 where $\tilde{\nu}_{n}(t)=n^{-1}\sum_{k=1}^{n}\tilde{B}_{k\Delta}t(X_{k\Delta})$.
By Lemma \ref{lem_majoration_AB-1}, we get that\[
\E{F_{k\Delta}^{2}}\lesssim\Delta+\sigma_{0}^{2}\Delta^{1-\beta/2}+\E{B_{k\Delta}^{2}\units{\Omega_{X,k}\cap\Omega_{N,k}^{c}}}+\left(\E{B_{k\Delta}\units{\Omega_{X,k}\cap\Omega_{N,k}}\left|\rond{F}_{k\Delta}\right.}\right)^{2}.\]
 By \eqref{eq:majoration_exp_Zk}, $\mathbb{P}\left(\left|Z_{k\Delta}\right|\geq\Delta^{1/2}\ln(n)\right)\lesssim n^{-1}$,
then $\mathbb{P}\left(\left|B_{k\Delta}\right|\geq\ln^{2}(n)\right)\lesssim n^{-1}$
and then: \begin{align*}
\E{B_{k\Delta}^{2}\units{\Omega_{X,k}\cap\Omega_{N,k}^{c}}} & \lesssim\E{B_{k\Delta}^{2}\units{\left|B_{k\Delta}\right|\geq\ln^{2}(n)}}+\E{\ln^{2}(n)\units{\Omega_{X,k}\cap\Omega_{N,k}^{c}}}\\
 & \lesssim n^{-1}+\ln^{2}(n)\mathbb{P}\left(\Omega_{X,k}\cap\Omega_{N,k}^{c}\right)\lesssim\ln^{2}(n)\Delta^{2-\beta/2}.\end{align*}
 As the random variables $B_{k\Delta}$ are centred:~ \begin{align*}
\left(\E{B_{k\Delta}\units{\Omega_{X,k}\cap\Omega_{N,k}}}\right)^{2} & =\left(\E{B_{k\Delta}\units{\Omega_{X,k}^{c}\cup\Omega_{N,k}^{c}}}\right)^{2}\\
 & \lesssim\E{B_{k\Delta}^{2}}\left(\mathbb{P}\left(\Omega_{X,k}^{c}\right)+\mathbb{P}\left(\Omega_{N,k}^{c}\right)\right)\\
 & \lesssim\Delta^{1-\beta/2}.\end{align*}
 Then $\E{F_{k\Delta}^{2}}\lesssim\Delta^{1-\beta/2}.$ We have that
\[
\E{\sup_{t\in S_{m}}\tilde{\nu}_{n}^{2}(t)}\leq\sum_{\lambda\in\Lambda}\E{\tilde{\nu}_{n}^{2}\left(\varphi_{\lambda}\right)}\leq\frac{\E{\tilde{B}_{k\Delta}^{2}}}{n}\lesssim\frac{\sigma_{0}^{4}}{n}.\]
 where$\left(\varphi_{\lambda}\right)_{1\leq\lambda\leq D_{m}}$ is
the orthonormal basis of $S_{m}$ for the $\norm{.}_{\pi}$-norm.

\subsection{Proof of Theorem \ref{thme_sigma_adaptatif}}

We apply the Berbee's coupling Lemma to the random exponentially $\beta$-mixing
vectors $(\tilde{B}_{k\Delta},X_{k\Delta})$. For any $a\in\{0,1\}$,
we can construct random variables \[
V_{k,a}^{*}=\frac{1}{q_{n}}\sum_{l=1}^{q_{n}}\tilde{B}_{(2(k-1)+a)q_{n}+l)\Delta}^{*}t(X_{(2(k-1+a)q_{n}+l)\Delta}^{*})\]
 independent and of same law as \[
V_{k,a}=\frac{1}{q_{n}}\sum_{l=1}^{q_{n}}\tilde{B}_{(2(k-1)+a)q_{n}+l)\Delta}t(X_{(2(k-1+a)q_{n}+l)\Delta}).\]
 Let us set $\tilde{\Omega}^{*}=\left\{ \omega,\:\forall a,\:\forall k,\: V_{k,a}=V_{k,a}^{*}\right\} $,
$\mathbb{P}\left(\tilde{\Omega}^{*c}\right)\lesssim n^{-4}$. Let
us consider the set $\Omega_{Z}=\left\{ \omega,\:\forall k\:\left|Z_{k\Delta}\right|\leq4\sigma_{0}\ln(n)\Delta^{1/2}\right\} $
on which the random variables $\tilde{B}_{k\Delta}$ are bounded.
According to inequality \eqref{eq:majoration_exp_Zk}, $\mathbb{P}\left(\Omega_{Z}^{c}\right)\lesssim n^{-4}$.

Let us set $\tilde{\rond{O}}=\Omega_{n}\cap\Omega_{Z}\cap\tilde{\Omega}^{*}$.
We bound the risk on $\tilde{\rond{O}}^{c}$ in the same way as on
$\Omega_{n}^{c}$. Let us set \[
\tilde{\nu}_{n}^{*}(t)=\tilde{\nu}_{n,0}^{*}(t)+\tilde{\nu}_{n,1}^{*}(t)\quad\textrm{with}\quad\tilde{\nu}_{n,a}^{*}(t)=\frac{1}{p_{n}}\sum_{k=1}^{p_{n}}V_{k,a}^{*}-\E{V_{k,a}^{*}}.\]
 For any $m\in\rond{M}_{n}$: \begin{align*}
\E{\norm{\hat{\sigma}_{\hat{m}}-\sigma_{A}}_{n}^{2}\units{\tilde{\rond{O}}}} & \leq3\norm{\sigma_{m}^{2}-\sigma^{2}}_{\pi}^{2}+12\E{F_{k\Delta}^{2}}+12\left(\E{\tilde{B}_{k\Delta}^{*}\units{\tilde{\rond{O}}}}\right)^{2}\\
 & +2pen(m)-2\widetilde{pen}(\hat{m})+\E{\sup_{t\in\rond{B}_{m,\hat{m}}}\left(\tilde{\nu}_{n}^{*}(t)\right)^{2}}.\end{align*}
 Let us introduce the function $\tilde{p}(m,m')=(\widetilde{pen}(m)+\widetilde{pen}(m'))/12.$
We have that \[
\left[\left(\sup_{t\in\rond{B}_{m,\hat{m}}}\tilde{\nu}_{n}^{*2}(t)-\tilde{p}(m,\hat{m})\right)\units{\tilde{\rond{O}}}\right]_{+}\leq\sum_{m'\in\rond{M}_{n}}\left[\left(\sup_{t\in\rond{B}_{m}}\tilde{\nu}_{n}^{*2}(t)-\tilde{p}(m,m')\right)\units{\tilde{\rond{O}}}\right]_{+}.\]
 On $\tilde{\rond{O}}$, for any $a$, the random variables $(V_{k,a}^{*})$
are independent, centred and bounded. We have that $\left|V_{k,a}^{*}\right|\leq\tilde{M}=\sigma_{0}^{2}\ln^{2}(n)D^{1/2}$,
$\E{\left(V_{k,a}^{*}\right)^{2}}\leq\tilde{V}=\sigma_{0}^{4}/q_{n}$
and \[
\E{\sup_{t\in\rond{B}_{m,m'}}\tilde{\nu}_{n,a}^{*2}(t)}\leq\tilde{H}=\sigma_{0}^{4}\frac{D}{n}.\]
 By the Talagrand's inequality, we deduce:~ \begin{align*}
R_{3}:= & \E{\left[\sup_{t\in\rond{B}_{m,m'}}\left(\tilde{\nu}_{n}^{*}(t)\right)^{2}-12\tilde{p}(m,m')\right]_{+}}\\
 & \lesssim\frac{1}{n}\exp\left(-c\frac{p_{n}q_{n}}{n}D\right)+\frac{\ln^{4}(n)D}{p_{n}^{2}}\exp\left(-c'\frac{p_{n}D^{1/2}}{\ln(n)D^{1/2}}\right)\\
 & \lesssim\frac{1}{n}\exp\left(-cD\right)+\frac{\ln^{6}(n)}{n^{2}\Delta^{2}}D\exp\left(-c'\frac{\sqrt{n\Delta}}{\ln^{3}(n)}\right).\end{align*}
 As $\rond{D_{n}^{2}}\leq n\Delta$ and $\ln^{3}(n)\ll n\Delta$ ,
we find:~ \begin{align*}
\E{\left[\sup_{t\in\rond{B}_{m,\hat{m}}}\nu_{n}^{2}(t)-pen(m)\right]_{+}} & \lesssim\frac{1}{n}\sum_{m'}e^{-cD_{m,m'}}+\frac{\ln^{4}(n)}{n\Delta}\exp-c'\frac{n\Delta}{\ln^{2}(n)}\\
 & \lesssim\frac{1}{n}.\end{align*}

\bibliographystyle{style}
\bibliography{biblio2}

\begin{thebibliography}{19}
\expandafter\ifx\csname natexlab\endcsname\relax\def\natexlab#1{#1}\fi
\expandafter\ifx\csname url\endcsname\relax
  \def\url#1{{\tt #1}}\fi
\expandafter\ifx\csname urlprefix\endcsname\relax\def\urlprefix{URL }\fi

\bibitem[{Applebaum(2004)}]{applebaum}
Applebaum, D. (2004) {\em L\'evy processes and stochastic calculus\/}, {\em
  Cambridge Studies in Advanced Mathematics\/}, volume~93.
\newblock Cambridge University Press, Cambridge.

\bibitem[{Arlot and Massart(2009)}]{arlotmassart2009}
Arlot, S. and Massart, P. (2009) Data-driven calbration of penalties for
  least-squares regression.
\newblock {\em Journal of Machine Learning Research\/}, 10 pp. 245--279.

\bibitem[{Birg{\'e} and Massart(1998)}]{birgemassart98}
Birg{\'e}, L. and Massart, P. (1998) Minimum contrast estimators on sieves:
  exponential bounds and rates of convergence.
\newblock {\em Bernoulli\/}, 4~(3) pp. 329--375.

\bibitem[{Comte and Rozenholc(2002)}]{comteroz2002}
Comte, F. and Rozenholc, Y. (2002) Adaptive estimation of mean and volatility
  functions in (auto-)regressive models.
\newblock {\em Stochastic Process. Appl.\/}, 97~(1) pp. 111--145.

\bibitem[{Comte and Rozenholc(2004)}]{comteroz2004}
Comte, F. and Rozenholc, Y. (2004) A new algorithm for fixed design regression
  and denoising.
\newblock {\em Ann. Inst. Statist. Math.\/}, 56~(3) pp. 449--473.

\bibitem[{Comte {\em et~al.\/}(2007)Comte, Genon-Catalot and
  Rozenholc}]{comtegenon2007}
Comte, F., Genon-Catalot, V. and Rozenholc, Y. (2007) Penalized nonparametric
  mean square estimation of the coefficients of diffusion processes.
\newblock {\em Bernoulli\/}, 13~(2) pp. 514--543.

\bibitem[{Comte and Merlev{\`e}de(2002)}]{comtemerlevede2002}
Comte, F. and Merlev{\`e}de, F. (2002) Adaptive estimation of the stationary
  density of discrete and continuous time mixing processes.
\newblock {\em ESAIM Probab. Statist.\/}, 6 pp. 211--238 (electronic).
\newblock New directions in time series analysis (Luminy, 2001).

\bibitem[{Dellacherie and Meyer(1980)}]{dellacheriemeyer}
Dellacherie, C. and Meyer, P.A. (1980) {\em Probabilit\'es et potentiel.
  {C}hapitres {V} \`a {VIII}\/}, {\em Actualit\'es Scientifiques et
  Industrielles [Current Scientific and Industrial Topics]\/}, volume 1385.
\newblock Hermann, Paris, revised edition.
\newblock Th{\'e}orie des martingales. [Martingale theory].

\bibitem[{DeVore and Lorentz(1993)}]{devorelorentz}
DeVore, R.A. and Lorentz, G.G. (1993) {\em Constructive approximation\/}, {\em
  Grundlehren der Mathematischen Wissenschaften [Fundamental Principles of
  Mathematical Sciences]\/}, volume 303.
\newblock Springer-Verlag, Berlin.

\bibitem[{Hanif {\em et~al.\/}(2012)Hanif, Wang and Lin}]{muhammad_wang_lin}
Hanif, M., Wang, H. and Lin, Z. (2012) Reweighted {N}adaraya-{W}atson
  estimation of jump-diffusion models.
\newblock {\em Sci. China Math.\/}, 55~(5) pp. 1005--1016.
\newblock \urlprefix\url{http://dx.doi.org/10.1007/s11425-011-4340-4}.

\bibitem[{Hoffmann(1999)}]{hofmann99}
Hoffmann, M. (1999) Adaptive estimation in diffusion processes.
\newblock {\em Stochastic Process. Appl.\/}, 79~(1) pp. 135--163.

\bibitem[{Kappus(2012)}]{kappus_these}
Kappus, J. (2012) {\em Nonparametric adaptive estimation for discretely
  observed L\'evy processes\/}.
\newblock Ph.D. thesis, Humboldt-Universit\"at zu Berlin.

\bibitem[{Mai(2012)}]{mai2012}
Mai, H. (2012) Efficient maximum likelihood estimation for lévy-driven
  ornstein-uhlenbeck processes.

\bibitem[{Mancini and Ren{\`o}(2011)}]{mancinireno2011}
Mancini, C. and Ren{\`o}, R. (2011) Threshold estimation of {M}arkov models
  with jumps and interest rate modeling.
\newblock {\em J. Econometrics\/}, 160~(1) pp. 77--92.

\bibitem[{Masuda(2007)}]{masuda2007}
Masuda, H. (2007) Ergodicity and exponential {$\beta$}-mixing bounds for
  multidimensional diffusions with jumps.
\newblock {\em Stochastic Process. Appl.\/}, 117~(1) pp. 35--56.

\bibitem[{Meyer(1990)}]{meyer}
Meyer, Y. (1990) {\em Ondelettes et op\'erateurs. {I}\/}.
\newblock Actualit\'es Math\'ematiques. [Current Mathematical Topics]. Hermann,
  Paris.
\newblock Ondelettes. [Wavelets].

\bibitem[{Rubenthaler(2010)}]{Rubenthaler_HDR}
Rubenthaler, S. (2010) {\em Probabilit{\'e}s : aspects th{\'e}oriques et
  applications en filtrage non lin{\'e}aire, syst{\`e}mes de particules et
  processus stochastiques.\/}.
\newblock Habilitation {\`a} diriger des recherches, Universit\'e de
  Nice-Sophia Antipolis, France.

\bibitem[{Shimizu(2008)}]{shimizu2008}
Shimizu, Y. (2008) Some remarks on estimation of diffusion coefficients for
  jump-diffusions from finite samples.
\newblock {\em Bull. Inform. Cybernet.\/}, 40 pp. 51--60.

\bibitem[{Viennet(1997)}]{viennet97}
Viennet, G. (1997) Inequalities for absolutely regular sequences: application
  to density estimation.
\newblock {\em Probab. Theory Related Fields\/}.

\end{thebibliography}

\end{document}